\documentclass[12pt,twoside]{article}
\usepackage[mathscr]{eucal}
\usepackage{amssymb}
\usepackage{amscd}
\usepackage{amsmath,amsthm}
\usepackage{latexsym}
\usepackage{indentfirst}
\usepackage{yhmath}

\def\R{\mathbb{R}}
\def\C{\mathbb{C}}

\def\N{\mathbb{N}}

\def\Z{\mathbb{Z}}

\def\ZZ{\mathfrak{Z}}
\def\gg{\mathfrak{g}}
\def\hh{\mathfrak{h}}
\def\ll{\mathfrak{l}}
\def\pp{\mathfrak{p}}
\def\kk{\mathfrak{k}}
\def\nn{\mathfrak{n}}

\def\mm{\mathfrak{m}}
\def\aa{\mathfrak{a}}
\def\tt{\mathfrak{t}}

\def\dd{\mathfrak{d}}

\def\oo{\mathfrak{o}}
\def\ss{\mathfrak{s}}
\def\ii{\mathfrak{i}}

\def\a{\alpha}

\def\g{\gamma}
\def\G{\Gamma}
\def\d{\delta}
\def\D{\Delta}
\def\e{\varepsilon}
\def\f{\varphi}

\def\l{\lambda}
\def\o{\omega}

\def\r{\rho}
\def\s{\sigma}

\def\t{\tau}

\def\pr{\prime}
\def\sm{\setminus}
\def\sub{\subseteq}

\def\Lra{\Leftrightarrow}
\def\lra{\longrightarrow}

\def\ov{\overline}

\def\cH{{\cal H}}

\def\cD{{\cal D}}

\def\cS{{\cal S}}

\def\cP{{\cal P}}

\def\cU{{\cal U}}

\def\qqu{\qquad}
\def\qu{\quad}
\def\fa{\forall}

\def\mr{\mathrm}
\def\le{\left}
\def\ri{\right}
\def\map{\mapsto}

\def\tim{\times}
\def\otim{\otimes}
\def\beg{\begin}
\def\ind{\indent}
\def\noind{\noindent}

\def\lb{\linebreak}

\def\fr{\frac}

\def\en{enumerate}

\def\wp{\wideparen}

\def\GL{\mr{GL}}

\def\rank{\mr{rank}}

\def\Int{\mr{Int}}

\def\SO{\mr{SO}}
\def\SU{\mr{SU}}
\def\Spin{\mr{Spin}}

\def\Tr{\mr{Tr}}

\def\ad{\mr{ad}}

\def\span{\mr{span}}

\begin{document}
\date{\empty}
\title{Corners and fundamental corners for the groups $\Spin(n,1)$}
\author{Domagoj Kova\v cevi\' c and Hrvoje Kraljevi\' c, University of Zagreb
\thanks{The authors were supported by the QuantiXLie Centre of Excellence, a project cofinanced by the Croatian Government and European Union through the European Regional Development Fund - the Competitiveness and Cohesion Operational Programme (Grant KK.01.1.1.01.0004).}
}
\maketitle

\noind\underline{Abstract.} We study corners and fundamental corners of the irreducible representations of the groups $G=\Spin(n,1)$ that are not elementary, i.e. that are equivalent to subquotients of reducible nonunitary principal seres representations. For even $n$ we obtain results in a way analogous to the results in [10] for the groups $\SU(n,1).$ Especially, we again get a bijection between the nonelementary part ${\hat G}^0$ of the unitary dual $\hat G$ and the unitary dual $\hat K.$ In the case of odd $n$ we get a bijection between ${\hat G}^0$ and a true subset of $\hat K.$ 

\section{Introduction}

\ind{\bf 1. Elementary representations.} Let $G$ be a connected semisimple Lie group with finite center,  $\gg_0$ its Lie algebra, $K$ its maximal compact subgroup, and $\gg_0=\kk_0\oplus\pp_0$ the corresponding Cartan decomposition of $\gg_0.$ Let\lb $P=MAN$ be a minimal parabolic subgroup of $G;$ here the Lie algebra $\aa_0$ of the subgroup $A$ is a Cartan subspace of $\pp_0,$ i.e. a Lie subalgebra of $\gg_0$ which is maximal among those contained in $\pp_0,$ $M=K\cap P$ is the centralizer of $\aa_0$ in $K,$ its Lie algebra $\mm_0$ is the centralizer of $\aa_0$ in $\kk_0,$ $N=\exp(\nn_0),$ where $\nn_0$ is the sum of root subspaces $\gg_0^{\a}$ with respect to some choice $\D^+(\gg_0,\aa_0)$ of positive restricted roots of the pair $(\gg_0,\aa_0).$ Denote by $\D_P$ the modular function of the group $P.$ Then $\D_P(m)=1$ for every $m\in M,$ $\D_P(n)=1$ for every $n\in N$ and for $H\in\aa_0$ we have
$$
\D_P(\exp\,H)=\mr{e}^{\Tr(\ad\,H)|\nn_0}=\mr{e}^{2\d(H)},\qu\d=\fr{1}{2}\sum_{\a\in\D^+(\gg_0,\aa_0)}(\dim\,\gg_0^{\a})\a.
$$ 
Thus,
$$
\D_P(man)=\mr{e}^{2\d(\log\,a)},\qu m\in M,\,a\in A,\,n\in N,
$$
where $\log:A\lra\aa_0$ is the inverse map of the bijection $\exp|\aa_0:\aa_0\lra A.$ Let $\s$ be an irreducible unitary representation of the compact group $M$ on a finitedimensional unitary space $\cH_{\s}.$ Let $\aa$ be the complexification of $\aa_0.$ For $\nu\in\aa^*$ let $a\map a^{\nu}$ be the onedimensional representation of the Abelian group $A$ defined by
$$
a^{\nu}=\mr{e}^{\nu(\log\,a)},\qu a\in A.
$$
Define the representation $\s\otim\nu$ of the group $P=MAN$ on the space $\cH_{\s}$ by
$$
(\s\otim\nu)(man)=a^{\nu}\s(m),\qu m\in M,\,a\in A,\,n\in N.
$$
Let $\pi^{\s,\nu}$ be the representation of $G$ induced by the representation $\s\otim\nu.$ The space of the representation $\pi^{\s,\nu}$ is the Hilbert space $\cH^{\s,\nu}$ of all (classes of) Haar$-$measurable functions $f:G\lra\cH_{\s}$ such that
$$
f(px)=\sqrt{\D_P(p)}(\s\otim\nu)(p)f(x)\qu\fa p\in P,\,\,\fa x\in G,
$$
and such that
$$
\int_K\|f(k)\|_{\cH_{\s}}^2\mr{d}\mu(k)<+\infty,
$$
where $\mu$ is the normed Haar measure on $K$ and $\|\,\cdot\,\|_{\cH_{\s}}$ is the norm on the unitary space $\cH_{\s}.$ The representation $\pi^{\s,\nu}$ is given by the right action of $G:$
$$
\le[\pi^{\s,\nu}(x)f\ri](y)=f(yx),\qu f\in\cH^{\s,\nu},\,x,y\in G.
$$
The representations $\pi^{\s,\nu},$ $\s\in\hat{M},$ $\nu\in\aa^*,$ are called {\bf elementary representations} of $G.$\\
\ind Since $\D_P(man)=a^{2\d},$ the condition $f(px)=\sqrt{\D_P(p)}(\s\otim\nu)(p)f(x)$ can be written as
$$
f(manx)=a^{\nu+\d}\s(m)f(x),\qu m\in M,\,a\in A,\,n\in N,\,x\in G.
$$
\ind From classical results of Harish$-$Chandra we know that all elementary representations are admissible and of finite length and that every completely irreducible admissible representation of $G$ on a Banach space is infinitesimally equivalent to an irreducible subquotient of an elementary representation. Infinitesimal equivalence of completely irreducible admissible representations means algebraic equivalence of the corresponding $(\gg,K)-$modules. We will denote by $\wp{G}$ the set of all infinitesimal equivalence classes of completely irreducible admissible representations of $G$ on Banach spaces. $\wp{G}^e$ will denote the set of infinitesimal equivalence classes of irreducible elementary representations and $\wp{G}^0=\wp{G}\sm\wp{G}^e$ the set of infinitesimal equivalence classes of irreducible suquotients of reducible elementary representations. It is also due to Harish$-$Chandra that every irreducible unitary representation is admissible and that infinitesimal equivalence between such representations is equivalent to their unitary equivalence. Thus the unitary dual $\hat{G}$ of $G$ can be regarded as a subset of $\wp{G}.$ We denote $\hat{G}^e=\hat{G}\cap\wp{G}^e$ and $\hat{G}^0=\hat{G}\cap\wp{G}^0=\hat{G}\sm\hat{G}^e.$\\

\ind{\bf 2. Infinitesimal characters.} For a finitedimensional complex Lie algebra $\gg$ we denote by $\cU(\gg)$ the universal enveloping algebra of $\gg$ and by $\ZZ(\gg)$ the center of $\cU(\gg).$ Any unital homomorphism $\chi:\ZZ(\gg)\lra\C$ is called {\bf infinitesimal character} of $\gg.$ We denote by $\hat{\ZZ}(\gg)$ the set of all infinitesimal characters of $\gg.$ If $\pi$ is a representation of $\gg$ on a vector space $V$ we say that $\chi\in\hat{\ZZ}(\gg)$ is the infinitesimal character of the representation $\pi$ (or of the corresponding $\cU(\gg)-$module $V)$ if
$$
\pi(z)v=\chi(z)v\qqu\fa z\in\ZZ(\gg),\,\,\fa v\in V.
$$
Let now $\gg$ be semisimple and let $\hh$ be its Cartan subalgebra. Denote by $\D=\D(\gg,\hh)\sub\hh^*$ the root system of the pair $(\gg,\hh),$ by $W=W(\gg,\hh)$ its Weyl group, by $\D^+$ a choice of positive roots in $\D,$ by $\gg^{\a}$ the root subspace of $\gg$ for a root $\a\in\D,$ and
$$
\nn=\sum_{\a\in\D^+}\dotplus\,\gg^{\a}\qqu\mr{and}\qqu\ov{\nn}=\sum_{\a\in\D^+}\dotplus\,\gg^{-\a}.
$$
Then we have direct sum decomposition
$$
\cU(\gg)=\cU(\hh)\,\dotplus\,(\nn\,\cU(\gg)+\cU(\gg)\,\ov{\nn}).
$$
Denote by $\eta:\cU(\gg)\lra\cU(\hh)$ the corresponding projection. By a result of Harish$-$Chandra the restriction $\eta|\ZZ(\gg)$ is an injective homomorphism of $\ZZ(\gg)$ into the algebra $\cU(\hh).$ Since the Lie algebra $\hh$ is Abelian, the algebra $\cU(\hh)$ identifies with the symmetric algebra $\cS(\hh)$ over $\hh,$ thus with the polynomial algebra $\cP(\hh^*)$ over the dual space $\hh^*$ of $\hh.$ Therefore $\eta|\ZZ(\gg)$ is a monomorphism of $\ZZ(\gg)$ into $\cP(\hh^*).$ This monomorphism depends on the choice of $\D^+.$ This dependence is repared by the automorphism $\g=\g_{\D^+}$ of the algebra $U(\hh)=\cP(\hh^*)$ defined by
$$
(\g(u))(\l)=u(\l-\r),\,\,\l\in\hh^*,\,\,u\in\cU(\hh)=\cP(\hh^*),\,\,\mr{where}\,\,\r=\r_{\D^+}=\fr{1}{2}\sum_{\a\in\D^+}\a.
$$ 
Now the restriction $\o=(\g\circ\eta)|\ZZ(\gg)$ is independent on the choice of $\D^+$ and is a unital isomorphism of the algebra $\ZZ(\gg)$ onto the algebra $\cP(\hh^*)^W$ of polynomial functions on $\hh^*$ invariant under the Weyl group $W=W(\gg,\hh)$ of the root system $\D=\D(\gg,\hh).$ $\o$ is called the {\bf Harish$-$Chandra isomorphism.} By evaluation at the points of $\hh^*$ one obtains all infinitesimal characters: for $\l\in\hh^*$ we define infinitesimal character $\chi_{\l}\in\hat{\ZZ}(\gg)$ by
$$
\chi_{\l}(z)=(\o(z))(\l)=(\eta(z))(\l-\r),\qqu z\in\ZZ(\gg).
$$ 
Then $\l\map\chi_{\l}$ is a surjection of $\hh^*$ onto $\hat{\ZZ}(\gg)$ and for $\l,\mu\in\hh^*$ we have $\chi_{\l}=\chi_{\mu}$ if and only if $\mu=w\l$ for some $w\in W.$\\
\ind A choice of an ordered basis $(H_1,\ldots,H_{\ell})$ of $\hh$ identifies the dual space $\hh^*$ with $\C^{\ell}:$ $\l\in\hh^*$ identifies with the $\ell-$tuple $(\l(H_1),\ldots,\l(H_{\ell}))\in\C^{\ell}.$ Now, if $\hh^{\pr}$ is another Cartan subalgebra of $\gg,$ then there exists an inner automorphism $\f$ of $\gg$ such that $\hh^{\pr}=\f(\hh).$ $\f$ carries $(H_1,\ldots,H_{\ell})$ to a basis $(H_1^{\pr},\ldots,H_{\ell}^{\pr})$ of $\hh^{\pr}$ which we use for the identification of ${\hh^{\pr}}^*$ with $\C^{\ell}.$ If an $\ell-$tuple $(c_1,\ldots,c_{\ell})\in\C^{\ell}$ corresponds to $\l\in\hh^*$ and to $\l^{\pr}\in{\hh^{\pr}}^*$ then the corresponding infinitesimal characters are the same: $\chi_{\l}=\chi_{\l^{\pr}}.$\\

\ind We return now to the notations of {\bf 1.} If $\ll_0$ is any real Lie algebra (or its subspace) we will denote by $\ll$ its complexification. It is well known that an elementary representation has infinitesimal character. We are going to write down the formula for the infinitesimal character of the elementary representation $\pi^{\s,\nu},$ $\s\in\hat{M},$ $\nu\in\aa^*.$ Let $\dd_0$ be a Cartan subalgebra of the reductive Lie subalgebra $\mm_0.$ Denote by $\D_{\mm}=\D(\mm,\dd)\sub\dd^*$ the root system of the pair $(\mm,\dd).$ Choose a subset $\D_{\mm}^+$ of positive roots in $\D_{\mm}$ and set
$$
\d_{\mm}=\r_{\D_{\mm}^+}=\fr{1}{2}\sum_{\a\in\D_{\mm}^+}\a.
$$
Denote by $\l_{\s}\in\dd^*$ the highest weight of the representation $\s$ with respect to $\D_{\mm}^+.$ Now, $\hh_0=\dd_0\dotplus\aa_0$ is a Cartan subalgebra of $\gg_0$ and its complexification $\hh=\dd\dotplus\aa$ is a Cartan subalgebra of $\gg.$ Then the infinitesimal character of the elementary representation $\pi^{\s,\nu}$ is $\chi_{\Lambda(\s,\nu)},$ where $\Lambda(\s,\nu)\in\hh^*$ is given by
$$
\Lambda(\s,\nu)|\dd=\l_{\s}+\d_{\mm}\qqu\mr{and}\qqu\Lambda(\s,\nu)|\aa=\nu.
$$

\ind{\bf 3. Corners and fundamental corners} Suppose now that the rank of $\gg$ is equal to the rank of $\kk.$ Choose a Cartan subalgebra $\tt_0$ of $\kk_0.$ It is then also Cartan subalgebra of $\gg_0$ and the complexification $\tt$ is Cartan subalgebra of the complexifications $\kk$ and $\gg.$ Let $\D_K=\D(\kk,\tt)\sub\D=\D(\gg,\tt)$ be the root systems of the pairs $(\kk,\tt)$ and $(\gg,\tt)$ and $W_K=W(\kk,\tt)\sub W=W(\gg,\tt)$ the corresponding Weyl groups. Choose positive roots $\D_K^+$ in $\D_K$ and let $C$ be the corresponding $W_K-$Weyl chamber in $\tt_{\R}^*=i\tt_0^*.$ Denote by $\cD$ the set of all $W-$Weyl chambers in $i\tt_0^*$ contained in $C.$ For $D\in\cD$ we denote by $\D^D$ the corresponding positive roots in $\D$ and let $\D_P^D$ be the noncompact roots in $\D^D,$ i.e. $\D_P^D=\D^D\sm\D_K^+.$ Set
$$
\r_K=\fr{1}{2}\sum_{\a\in\D_K^+}\a\qqu\mr{and}\qqu\r_P^D=\fr{1}{2}\sum_{\a\in\D_P^D}\a.
$$
\ind Recall some definitions from [10]. For a representation $\pi$ of $G$ and for $q\in\hat{K}$ we denote by $(\pi:q)$ the multiplicity of $q$ in $\pi|K.$ The $K-${\bf spectrum} $\G(\pi)$ of a representation $\pi$ of $G$ is defined by
$$
\G(\pi)=\{q\in\hat{K};\,\,(\pi:q)>0\}.
$$
We identify $q\in\hat{K}$ with its maximal weight in $i\tt_0^*$ with respect to $\D_K^+.$ For $q\in\G(\pi)$ and for $D\in\cD$ we say: 
\beg{\en}
\item[$(i)$] $q$ is a $D-${\bf corner} for $\pi$ if $q-\a\not\in\G(\pi)$ $\fa\a\in\D_P^D;$
\item[$(ii)$] $q$ is a $D-${\bf fundamental corner} for $\pi$ if it is a $D-$corner for $\pi$ and $\chi_{q+\r_K-\r_P^D}$ is the infinitesimal character of $\pi;$
\item[$(iii)$] $q$ is {\bf fundamental corner} for $\pi$ if it is a $D-$fundamental corner for $\pi$ for some $D\in\cD.$
\end{\en}
\ind In [10] for the case of the groups $G=SU(n,1)$ and $K=U(n)$ the following results were proved:
\beg{\en}
\item[{\bf 1.}] Elementary representation $\pi^{\s,\nu}$ is reducible if and only if there exist $q\in\G(\pi^{\s,\nu})$ and $D\in\cD$ such that $\chi_{q+\r_K-\r_P^D}$ is the infinitesimal character of $\pi^{\s,\nu},$ i.e. if and only if $\Lambda(\s,\nu)=w(q+\r_K-\r_P^D)$ for some $w\in W.$
\item[{\bf 2.}] Every $\pi\in\wp{G}^0$ has either one or two fundamental corners.
\item[{\bf 3.}] $\hat{G}^0=\{\pi\in\wp{G}^0;\,\,\pi\,\,\mr{has}\,\,\mr{exactly}\,\,\mr{one}\,\,\mr{fundamental}\,\,\mr{corner}\}.$
\item[{\bf 4.}] For $\pi\in\hat{G}^0$ denote by $q(\pi)$ the unique fundamental corner of $\pi.$ Then $\pi\map q(\pi)$ is a bijection of $\hat{G}^0$ onto $\hat{K}.$
\end{\en}
\ind In this paper we investigate the analogous notions and results for the groups $\Spin(n,1).$

\section{The groups $\Spin(n,1)$}

\ind In the rest of the paper $G=\Spin(n,1),$ $n\geq3,$ is the connected and simply connected real Lie group with simple real Lie algebra
$$
\gg_0=\ss\oo(n,1)=\{A\in\gg\ll(n+1,\R);\,\,A^t=-\G A\G\},\qqu\G=\le[\beg{array}{cc}I_n&0\\0&-1\end{array}\ri],
$$
i.e.
$$
\gg_0=\le\{\le[\beg{array}{cc}B&a\\a^t&0\end{array}\ri];\,\,B\in\ss\oo(n),\,\,a\in M_{n,1}(\R)\ri\}.
$$
Here and in the rest of the paper we use the usual notation:
\beg{\en}
\item[$\bullet$] For $n,m\in\N$ $M_{m,n}(K)$ is the vector space of $m\tim n$ matrices over a field $K.$
\item[$\bullet$] $\gg\ll(n,K)$ is $M_{n,n}(K),$ considered as a Lie algebra with commutator $[A,B]=AB-BA.$
\item[$\bullet$] $\GL(n,K)$ is the group of invertible matrices in $M_{n,n}(K).$
\item[$\bullet$] $A^t$ is the transpose of a matrix $A.$
\item[$\bullet$] $\ss\oo(n,K)=\{B\in\gg\ll(n,K);\,\,B^t=-B\}.$
\item[$\bullet$] $\ss\oo(n)=\ss\oo(n,\R).$
\item[$\bullet$] $\SO(n)=\{A\in\GL(n,\R);\,\,A^{-1}=A^t,\,\,\det\,A=1\}.$
\end{\en}

\newpage

\ind For the group $G=\Spin(n,1)$ we choose Cartan decomposition $\gg_0=\kk_0\oplus\pp_0$ as follows
$$
\kk_0=\le\{\le[\beg{array}{cc}B&0\\0&0\end{array}\ri];\,\,B\in\ss\oo(n)\ri\},\qqu\pp_0=\le\{\le[\beg{array}{cc}0&a\\a^t&0\end{array}\ri];\,\,a\in M_{n,1}(\R)\ri\}.
$$
The complexifications are:
$$
\gg=\ss\oo(n,1,\C)=\{A\in\gg\ll(n+1,\C);\,\,A^t=-\G A\G\},
$$
i.e.
$$
\gg=\le\{\le[\beg{array}{cc}B&a\\a^t&0\end{array}\ri];\,\,B\in\ss\oo(n,\C),\,\,a\in M_{n,1}(\C)\ri\},
$$
$$
\kk=\le\{\le[\beg{array}{cc}B&0\\0&0\end{array}\ri];\,\,B\in\ss\oo(n,\C)\ri\},\qqu\pp=\le\{\le[\beg{array}{cc}0&a\\a^t&0\end{array}\ri];\,\,a\in M_{n,1}(\C)\ri\}.
$$
$\Spin(n,1)$ is double cover of the identity component $\SO_0(n,1)$ of the Lie group
$$
\SO(n,1)=\{A\in\GL(n+1,\R);\,\,A^{-1}=\G A^t\G,\,\,\det\,A=1\}.
$$
The analytic subgroup $K\subset G$ whose Lie algebra is $\kk_0$ is a maximal compact subgroup of $G$ isomorphic with the double cover $\Spin(n)$ of the group $\SO(n).$\\
\ind Now we choose Cartan subalgebras. $E_{p,q}$ will denote the $(n+1)\tim(n+1)$ matrix with $(p,q)-$entry equal $1$ and all the other entries $0.$ Set
$$
I_{p,q}=E_{p,q}-E_{q,p},\qu1\leq p,q\leq n,\qu p\not=q,
$$
and
$$
B_p=E_{p,n+1}+E_{n+1,p},\qu1\leq p\leq n.
$$
Then $\{I_{p,q};\,\,1\leq q<p\leq n\}$ is a basis of the real Lie algebra $\kk_0$ and of its complexification $\kk$ and $\{B_p;\,\,1\leq p\leq n\}$ is a basis of the real subspace $\pp_0$ and of its complexification $\pp.$ Now $\tt_0=\span_{\R}\le\{I_{2p,2p-1};\,\,1\leq p\leq\fr{n}{2}\ri\}$ is a Cartan subalgebra of $\kk_0$ and its complexification $\tt=\span_{\C}\le\{I_{2p,2p-1};\,\,1\leq p\leq\fr{n}{2}\ri\}$ is a Cartan subalgebra of $\kk.$\\
\ind We consider now separately two cases: $n$ even and $n$ odd.

\newpage

\beg{center}
{\bf $n$ even, $n=2k$}
\end{center}

\ind In this case $\tt_0$ is also a Cartan subalgebra of $\gg_0$ and $\tt$ is a Cartan subalgebra of $\gg.$ Set
$$
H_p=-iI_{2p,2p-1},\qqu1\leq p\leq k.
$$
Dual space $\tt^*$ identifies with $\C^k$ as follows:
$$
\tt^*\ni\l=(\l(H_1),\ldots,\l(H_k))\in\C^k.
$$
Let $\{\a_1,\ldots,\a_k\}$ be the canonical basis of $\C^k=\tt^*.$ The root system of the pair $(\gg,\tt)$ is
$$
\D=\D(\gg,\tt)=\{\pm\a_p\pm\a_q;\,\,1\leq p,q\leq k,\,\,p\not=q\}\cup\{\pm\a_p;\,\,1\leq p\leq k\}.
$$
The Weyl group $W$ of $\D$ consists of all permutations of the coordinates combined with multiplying some coordinates with $-1:$
$$
W=\Z_2^k\rtimes S_k=\{(\e,\s);\,\,\e\in\Z_2^k,\,\,\s\in S_k\},
$$
where $\Z_2$ is the multiplicative group $\{1,-1\}$ and $S_k$ is the group of permutations of $\{1,\ldots,k\}.$ $(\e,\s)\in W$ acts on $\tt^*=\C^k$ as follows:
$$
(\e,\s)(\l_1,\l_2,\ldots,\l_k)=(\e_1\l_{\s(1)},\e_2\l_{\s(2)},\ldots,\e_k\l_{\s(k)}).
$$
\ind The root system $\D_K$ of the pair $(\kk,\tt)$ is $\{\pm\a_p\pm\a_q;\,\,p\not=q\}.$ We choose positive roots in $\D_K:$
$$
\D_K^+=\{\a_p\pm\a_q;\,\,1\leq p<q\leq k\}.
$$
The corresponding Weyl chamber in $\R^k=i\tt_0^*$ is
$$
C=\{\l\in\R^k;\,\,(\l|\g_j)>0,\,1\leq j\leq k\}=\{\l\in\R^k;\,\l_1>\cdots>\l_{k-1}>|\l_k|>0\},
$$
and its closure is
$$
\ov{C}=\{\l\in\R^k;\,(\l|\g_j)\geq0,\,1\leq j\leq k\}=\{\l\in\R^k;\,\l_1\geq\cdots\geq\l_{k-1}\geq|\l_k|\}.
$$
The Weyl group $W_K$ of the root system $\D_K$ is the subgroup of $W$ consisting of all $(\e,\s)$ with even number of $\e_j=-1:$
$$
W_K=\{(\e,\s)\in W;\,\,\e_1\e_2\cdots\e_k=1\}\simeq\Z_2^{k-1}\rtimes S_k.
$$
We parametrize now the equivalence classes of irreducible finitedimensional representations of the Lie algebra $\kk$ (i.e. the unitary dual $\hat{K}$ of the group $K=\Spin(2k)\,)$ by identifying them with the corresponding highest weights. Thus
$$
\beg{array}{c}
\hat{K}=\le\{(m_1,\ldots,m_k)\in\Z^k\cup\le(\fr{1}{2}+\Z\ri)^k;\,\,m_1\geq m_2\geq\cdots\geq m_{k-1}\geq|m_k|\ri\}.
\end{array}
$$

\beg{center}
{\bf $n$ odd, $n=2k+1$}
\end{center}

\ind Now $\tt_0$ is not a Cartan subalgebra of $\gg_0.$ Set
$$
H=B_n=B_{2k+1}=E_{2k+1,2k+2}+E_{2k+2,2k+1},\qu\aa_0=\R H,\qu\hh_0=\tt_0\dotplus\aa_0.
$$
Then $\hh_0$ is a Cartan subalgebra of $\gg_0$ and all the other Cartan subalgebras of $\gg_0$ are $\Int(\gg_0)-$conjugated with $\hh_0.$ The ordered basis $(H_1,\ldots,H_k,H)$ of the complexification $\hh$ of $\hh_0$ is used for the identificaton of $\hh^*$ with $\C^{k+1}:$
$$
\hh^*\ni\l=(\l(H_1),\ldots,\l(H_k),\l(H))\in\C^{k+1}.
$$
$\tt^*$ identifies with $\C^k$ through ordered basis $(H_1,\ldots,H_k)$ of $\tt$ and $\aa^*$ identifies with $\C$ through $H:$
$$
\tt^*\ni\mu=(\mu(H_1),\ldots,\mu(H_k))\in\C^k,\qqu\aa^*\ni\nu=\nu(H)\in\C.
$$
Furthermore, $\tt^*$ and $\aa^*$ are identified with subspaces of $\hh^*$ as follows:
$$
\tt^*=\{\l\in\hh^*;\,\,\l|\aa=0\}=\{\l\in\C^{k+1};\,\,\l_{k+1}=0\},
$$
$$
\aa^*=\{\l\in\hh^*;\,\,\l|\tt=0\}=\{(0,\ldots,0,\nu);\,\,\nu\in\C\}.
$$
\ind Let $\{\a_1,\ldots,\a_{k+1}\}$ be the canonical basis of $\C^{k+1}.$ The root system\lb $\D=\D(\gg,\hh)$ of the pair $(\gg,\hh)$ is
$$
\D=\{\pm\a_p\pm\a_q;\,\,1\leq p,q\leq k+1,\,\,p\not=q\}.
$$
The Weyl group $W=W(\gg,\hh)$ consists of all permutations of coordinates combined with multiplying even number of coordinates with $-1:$
$$
W=\Z_2^k\rtimes S_{k+1}=\{(\e,\s);\,\,\e\in\Z_2^{k+1},\,\,\e_1\cdots\e_{k+1}=1,\,\,\s\in S_{k+1}\}.
$$
The root system $\D_K=\D(\kk,\tt)$ of the pair $(\kk,\tt)$ is
$$
\D_K=\{\pm\a_p\pm\a_q;\,\,1\leq p,q\leq k,\,\,p\not=q\}\cup\{\pm\a_p;\,\,1\leq p\leq k\}.
$$
Choose positive roots in $\D_K$ as follows:
$$
\D_K^+=\{\a_p\pm\a_q;\,\,1\leq p<q\leq k\}\cup\{\a_p;\,\,1\leq p\leq k\}.
$$
The corresponding Weyl chamber in $\R^k=i\tt_0^*$ is
$$
C=\{\l\in\R^k;\,(\l|\g_j)>0,\,1\leq j\leq k\}=\{\l\in\R^k;\,\l_1>\cdots>\l_k>0\}
$$
and its closure is
$$
\ov{C}=\{\l\in\R^k;\,(\l|\g_j)\geq0,\,1\leq j\leq k\}=\{\l\in\R^k;\,\l_1\geq\cdots\geq\l_k\geq0\}
$$
The dual $\hat{K}$ is again identified with the highest weights of ireducible representations. Thus:
$$
\beg{array}{c}
\hat{K}=\le\{q=(m_1,\ldots,m_k)\in\Z_+^k\cup\le(\fr{1}{2}+\Z_+\ri)^k;\,\,m_1\geq m_2\geq\cdots\geq m_k\ri\}.
\end{array}
$$

\beg{center}
{\bf Elementary representations of the groups $\Spin(n,1)$}
\end{center}

\ind Regardless the parity of $n$ we put
$$
H=B_n=E_{n,n+1}+E_{n+1,n},\qqu\aa_0=\R H.
$$
Then $\aa_0$ is maximal among Abelian subalgebras of $\gg_0$ contained in $\pp_0.$ As we already said, if $n$ is odd, $n=2k+1,$ then $\hh_0=\tt_0\dotplus\aa_0$ is a Cartan subalgebra of $\gg_0$ and all the other Cartan subalgebras are $\Int(\gg_0)-$conjugated to $\hh_0.$ If $n$ is even, $n=2k,$ set
$$
\hh_0=\span_{\R}\{iH_1,\ldots,iH_{k-1},H\}.
$$ 
It is a Cartan subalgebra of $\gg_0.$ In this case $\gg_0$ has two $\Int(\gg_0)-$conjugacy classes of Cartan subalgebras; $\hh_0$ and $\tt_0$ are their representatives. Their complexifications $\hh$ and $\tt$ are $\Int(\gg)-$conjugated. Explicitely, the matrix
$$
A=\le[\beg{array}{ccc}\fr{1}{\sqrt{2}}P_k&\fr{1}{\sqrt{2}}P_k&-ie_k\\-\fr{1}{\sqrt{2}}Q_k&\fr{1}{\sqrt{2}}I_k&0_k\\-\fr{i}{\sqrt{2}}e_k^t&\fr{i}{\sqrt{2}}e_k^t&0\end{array}\ri]\in\SO(2k,1,\C),
$$
where $P_k=I_k-E_{k,k}=\mr{diag}(1,\ldots,1,0),$ $Q_k=I_k-2E_{k,k}=\mr{diag}(1,\ldots,1,-1),$ $e_k\in M_{k,1}(\C)$ is given by $e_k^t=[0\cdots0\,1]$ and $0_k$ is the zero matrix in $M_{k,1}(\C),$ has the properties
$$
AH_jA^{-1}=H_j,\qu1\leq j\leq k-1,\qu\mr{and}\qu AH_kA^{-1}=H;
$$
thus, $A\tt A^{-1}=\hh.$ As we mentioned before, this means that the parameters from $\C^k=\hh^*=\tt^*$ of the infinitesimal characters obtained through the two Harish$-$Chandra isomorphisms $\ZZ(\gg)\lra\cP(\hh^*)^W$ and $\ZZ(\gg)\lra\cP(\tt^*)^W$ coincide if the identifications of $\hh^*$ and $\tt^*$ with $\C^k$ are done throught the two ordered bases $(H_1,\ldots,H_{k-1},H)$ of $\hh$ and $(H_1,\ldots,H_{k-1},H_k)$ of $\tt.$\\

\ind For both cases, $n$ even and $n$ odd, $\mm_0$ (the centralizer of $\aa_0$ in $\kk_0)$ is the subalgebra of all matrices in $\gg_0$ with the last two rows and columns $0.$ The subgroup $M$ is isomorphic to $\Spin(n-1).$ A Cartan subalgebra of $\mm_0$ is
$$
\dd_0=\tt_0\cap\mm_0=\span_{\R}\{iH_1,\ldots,iH_{k-1}\},\qqu k=\le[\fr{n}{2}\ri].
$$
The elements of $\hat{M}$ are identified with their highest weights. For $n$ even, $n=2k,$ we have
$$
\beg{array}{c}
\hat{M}=\le\{(n_1,\ldots,n_{k-1})\in\Z_+^{k-1}\cup\le(\fr{1}{2}+\Z_+\ri)^{k-1};\,\,n_1\geq n_2\geq\cdots\geq n_{k-1}\geq0\ri\}
\end{array}
$$
and for $n$ odd, $n=2k+1,$ we have
$$
\beg{array}{c}
\end{array}
\hat{M}=\le\{(n_1,\ldots,n_k)\in\Z^k\cup\le(\fr{1}{2}+\Z\ri)^k;\,\,n_1\geq n_2\geq\cdots\geq n_{k-1}\geq|n_k|\ri\}.
$$
The branching rules for the restriction of representations of $K$ to the subgroup $M$ are the following:\\
\ind If $n$ is even, $n=2k,$ we have
$$
(m_1,\ldots,m_k)|M=\bigoplus_{(n_1,\ldots,n_{k-1})\prec(m_1,\ldots,m_k)}(n_1,\ldots,n_{k-1});
$$
here the symbol $(n_1,\ldots,n_{k-1})\prec(m_1,\ldots,m_k)$ means that either all $m_i$ and $n_j$ are in $\Z$ or all of them are in $\fr{1}{2}+\Z$ and
$$
m_1\geq n_1\geq m_2\geq n_2\cdots\geq m_{k-1}\geq n_{k-1}\geq|m_k|.
$$
\ind If $n$ is odd, $n=2k+1,$ we have
$$
(m_1,\ldots,m_k)|M=\bigoplus_{(n_1,\ldots,n_k)\prec(m_1,\ldots,m_k)}(n_1,\ldots,n_k);
$$
now the symbol $(n_1,\ldots,n_k)\prec(m_1,\ldots,m_k)$ means again that either all $m_i$ and $n_j$ are in $\Z$ or all of them are in $\fr{1}{2}+\Z$ and now that 
$$
m_1\geq n_1\geq m_2\geq n_2\cdots\geq m_{k-1}\geq n_{k-1}\geq m_k\geq|n_k|.
$$
\ind The restriction $\pi^{\s,\nu}|K$ is the representation of $K$ induced by the representation $\s$ of the subgroup $M,$ thus it does not depend on $\nu.$ By Frobenius Reciprocity Theorem the multiplicity of $q\in\hat{K}$ in $\pi^{\s,\nu}|K$ is equal to the multiplicity of $\s$ in $q|M.$ Thus
$$
\pi^{\s,\nu}|K=\bigoplus_{\s\prec(m_1,\ldots,m_k)}(m_1,\ldots,m_k).
$$
Hence, the multiplicity of every $q=(m_1,\ldots,m_k)\in\hat{K}$ in the elementary representation $\pi^{\s,\nu}$ is either $1$ or $0$ and the $K-$spectrum $\G(\pi^{\s,\nu})$ consists of all $q=(m_1,\ldots,m_k)\in\hat{K}\cap(n_1+\Z)^k$ such that
$$
\beg{array}{ll}m_1\geq n_1\geq m_2\geq n_2\geq\cdots\geq m_{k-1}\geq n_{k-1}\geq|m_k|&\,\,\mr{if}\,\,n=2k\\
&\\
m_1\geq n_1\geq m_2\geq n_2\geq\cdots\geq m_{k-1}\geq n_{k-1}\geq m_k\geq|n_k|&\,\,\mr{if}\,\,n=2k+1.
\end{array}
$$

\section{Representations of $\Spin(2k,1)$}

\ind In this section we first write down in our notation the known results on elementary representations and its irreducible subquotients for the groups $\Spin(2k,1)$ (see [1], [2], [3], [7], [8], [9], [11], [12]). For $\s=(n_1,\ldots,n_{k-1})$ in $\hat{M}\sub\R^{k-1}=i\dd_0^*$ and for $\nu\in\C=\aa^*$ the elementary representation $\pi^{\s,\nu}$ is irreducible if and only if either $\nu\not\in\fr{1}{2}+n_1+\Z$ or
$$
\beg{array}{c}
\nu\in\le\{\pm\le(n_{k-1}+\fr{1}{2}\ri),\pm\le(n_{k-2}+\fr{3}{2}\ri),\ldots,\pm\le(n_2+k-\fr{5}{2}\ri),\pm\le(n_1+k-\fr{3}{2}\ri)\ri\}.
\end{array}
$$
If $\pi^{\s,\nu}$ is reducible it has either two or three irreducible subquotients. If it has two, we will denote them by $\t^{\s,\nu}$ and $\o^{\s,\nu};$ an exception is the case of nonintegral $n_j$ and $\nu=0,$ when we denote them by $\o^{\s,0,\pm}.$ If $\pi^{\s,\nu}$ has three irreducible subquotients, we will denote them by $\t^{\s,\nu}$ and $\o^{\s,\nu,\pm}.$ Their $K-$spectra are as follows:
\beg{\en}
\item[$(a1)$] If $n_j\in\Z_+$ and $\nu\in\le\{\pm\fr{1}{2},\pm\fr{3}{2},\ldots,\pm\le(n_{k-1}-\fr{1}{2}\ri)\ri\}$ (this is possible only if $n_{k-1}\geq1)$ the representation $\pi^{\s,\nu}$ has three irreducible subquotients $\t^{\s,\nu}$ and $\o^{\s,\nu,\pm}.$ Their $K-$spectra consist of all $q=(m_1,\ldots,m_k)$ in $\hat{K}\cap\Z^k$ such that:
$$
\beg{array}{ll}
\G(\t^{\s,\nu}):&\,\,m_1\geq n_1\geq\cdots\geq m_{k-1}\geq n_{k-1},\,\,|m_k|\leq|\nu|-\fr{1}{2};\\
\G(\o^{\s,\nu,+}):&\,\,m_1\geq n_1\geq\cdots\geq m_{k-1}\geq n_{k-1}\geq m_k\geq|\nu|+\fr{1}{2};\\
\G(\o^{\s,\nu,-}):&\,\,m_1\geq n_1\geq\cdots\geq m_{k-1}\geq n_{k-1},\,\,\,-|\nu|-\fr{1}{2}\geq m_k\geq -n_{k-1}.
\end{array}
$$
\item[$(a2)$] If $n_j\in\le(\fr{1}{2}+\Z_+\ri)$ and $\nu\in\le\{\pm1,\ldots,\pm\le(n_{k-1}-\fr{1}{2}\ri)\ri\}$ (this is possible only if $n_{k-1}\geq\fr{3}{2})$ the representation $\pi^{\s,\nu}$ has three irreducible subquotients $\t^{\s,\nu}$ and $\o^{\s,\nu,\pm}.$ Their $K-$spectra consist of all $q=(m_1,\ldots,m_k)$ in $\hat{K}\cap\le(\fr{1}{2}+\Z\ri)^k$ such that:
$$
\beg{array}{ll}
\G(\t^{\s,\nu}):&\,\,m_1\geq n_1\geq\cdots\geq m_{k-1}\geq n_{k-1},\,\,|m_k|\leq|\nu|-\fr{1}{2};\\
\G(\o^{\s,\nu,+}):&\,\,m_1\geq n_1\geq\cdots\geq m_{k-1}\geq n_{k-1}\geq m_k\geq|\nu|+\fr{1}{2};\\
\G(\o^{\s,\nu,-}):&\,\,m_1\geq n_1\geq\cdots\geq m_{k-1}\geq n_{k-1},\,\,\,-|\nu|-\fr{1}{2}\geq m_k\geq -n_{k-1}.
\end{array}
$$
\item[$(a3)$] If $n_j\in\le(\fr{1}{2}+\Z_+\ri)$ and if $\nu=0$ the representation has two irreducible subquotients $\o^{\s,0,\pm};$ they are both subrepresentations since $\pi^{\s,0}$ is unitary. Their $K-$spectra consist of all $q=(m_1,\ldots,m_k)$ in $\hat{K}\cap\le(\fr{1}{2}+\Z\ri)^k$ such that:
$$
\beg{array}{ll}
\G(\o^{\s,0,+}):&\,\,m_1\geq n_1\geq\cdots\geq m_{k-1}\geq n_{k-1}\geq m_k\geq\fr{1}{2};\\
\G(\o^{\s,0,-}):&\,\,m_1\geq n_1\geq\cdots\geq m_{k-1}\geq n_{k-1},\,\,\,-\fr{1}{2}\geq m_k\geq -n_{k-1}.
\end{array}
$$
\item[$(bj)$] If $n_{j-1}>n_j$ for some $j\in\{2,\ldots,k-1\}$ and if
$$
\beg{array}{c}
\nu\in\le\{\pm\le(n_j+k-j+\fr{1}{2}\ri),\pm\le(n_j+k-j+\fr{3}{2}\ri),\ldots,\pm\le(n_{j-1}+k-j-\fr{1}{2}\ri)\ri\},
\end{array}
$$
then $\pi^{\s,\nu}$ has two irreducible subquotients $\t^{\s,\nu}$ and $\o^{\s,\nu}.$ Their $K-$spectra consist of all $q=(m_1,\ldots,m_k)\in\hat{K}\cap(n_1+\Z)^k$ such that:
$$
\beg{array}{ll}
\G(\t^{\s,\nu}):&\qu m_1\geq n_1\geq\cdots\geq m_{j-1}\geq n_{j-1},\\
&\qu|\nu|-k+j-\fr{1}{2}\geq m_j\geq n_j\geq\cdots\geq n_{k-1}\geq|m_k|;\\
\G(\o^{\s,\nu}):&\qu m_1\geq n_1\geq\cdots\geq n_{j-1}\geq m_j\geq|\nu|-k+j+\fr{1}{2},\\
&\qu n_j\geq m_{j+1}\geq\cdots\geq n_{k-1}\geq|m_k|.
\end{array}
$$
\item[$(c)$] If
$$
\beg{array}{c}
\nu\in\le\{\pm\le(n_1+k-\fr{1}{2}\ri),\pm\le(n_1+k+\fr{1}{2}\ri),\pm\le(n_1+k+\fr{3}{2}\ri),\ldots\ri\},
\end{array}
$$
then $\pi^{\s,\nu}$ has two irreducible subquotients: finitedimensional representation $\t^{\s,\nu}$ and infinitedimensional $\o^{\s,\nu}.$ Their $K-$spectra consist of all $q=(m_1,\ldots,m_k)\in\hat{K}\cap\le(n_1+\Z\ri)^k$ such that:
$$
\beg{array}{ll}
\G(\t^{\s,\nu}):&\,\,|\nu|-k+\fr{1}{2}\geq m_1\geq n_1\geq\cdots\geq m_{k-1}\geq n_{k-1}\geq|m_k|;\\
\G(\o^{\s,\nu}):&\,\,m_1\geq|\nu|-k+\fr{3}{2},\,\,n_1\geq m_2\geq\cdots\geq m_{k-1}\geq n_{k-1}\geq|m_k|.
\end{array}
$$
\end{\en}

\ind Irreducible elementary representation $\pi^{\s,\nu}$ is unitary if and only if either $\nu\in i\R$ (so called {\bf unitary principal series}) or $\nu\in\langle-\nu(\s),\nu(\s)\rangle,$ where
$$
\nu(\s)=\min\,\{\nu\geq0;\,\,\pi^{\s,\nu}\,\,\mr{is}\,\,\mr{reducible}\}
$$
(so caled {\bf complementary series}). Notice that for nonintegral $n_j$'s $\pi^{\s,0}$ is reducible, thus $\nu(\s)=0$ and the complementary series is empty. In the case of integral $n_j$'s we have the following possibilities:
\beg{\en}
\item[$(a)$] If $n_{k-1}\geq1,$ then $\nu(\s)=\fr{1}{2}.$ The reducible elementary representation $\pi^{\s,\fr{1}{2}}$ is of the type $(a1).$
\item[$(b)$] If $n_{k-1}=0$ and $n_1\geq1,$ let $j\in\{2,\ldots,k-1\}$ be such that\lb $n_{k-1}=\cdots=n_j=0<n_{j-1}.$ Then $\nu(\s)=k-j+\fr{1}{2}.$ The reducible elementary representation $\pi^{\s,k-j+\fr{1}{2}}$ is of the type $(bj).$
\item[$(c)$] If $\s$ is trivial, i.e. $n_1=\cdots=n_{k-1}=0,$ then $\nu(\s)=k-\fr{1}{2}.$ The reducible elementary representation $\pi^{\s,k-\fr{1}{2}}$ is of the type $(c).$
\end{\en}

\ind Among irreducible subquotients of reducible elementary representations the unitary ones are $\o^{\s,\nu,\pm},$ $\t^{\s,\nu(\s)}$ and $\o^{\s,\nu(\s)}.$\\

\ind Now we write down the infinitesimal characters. The dual space $\hh^*$ of the Cartan subalgebra $\hh=\dd\dotplus\aa$ is identified with $\C^k$ through the basis $(H_1,\ldots,H_{k-1},H).$ The infinitesimal character of the elementary representation $\pi^{\s,\nu}$ is $\chi_{\Lambda(\s,\nu)},$ where $\Lambda(\s,\nu)\in\hh^*$ is given by
$$
\beg{array}{c}
\Lambda(\s,\nu)|\dd=\l_{\s}+\d_{\mm}\qqu\mr{and}\qqu\Lambda(\s,\nu)|\aa=\nu,
\end{array}
$$
where $\l_{\s}\in\dd^*$ is the highest weight of the representation $\s$ and $\d_{\mm}$ is the halfsum of positive roots of the pair $(\mm,\dd).$ Using the earlier described identifications of $\dd^*=\C^{k-1}$ and $\aa^*=\C$ with subspaces of $\hh^*=\C^k$ we have\lb $\l_{\s}=(n_1,\ldots,n_{k-1},0),$ $\d_{\mm}=\le(k-\fr{3}{2},k-\fr{5}{2},\ldots,n_{k-1}+\fr{1}{2},0\ri),$ $\nu=(0,\ldots,0,\nu),$ hence
$$
\beg{array}{c}
\Lambda(\s,\nu)=\le(n_1+k-\fr{3}{2},n_2+k-\fr{5}{2},\ldots,n_{k-1}+\fr{1}{2},\nu\ri).
\end{array}
$$
As we pointed out, if $\tt^*$ is identified with $\C^k$ through the basis $(H_1,\ldots,H_k)$ od $\tt,$ the same parameters determine this infinitesimal character with respect to Harish$-$Chandra isomorphism $\ZZ(\gg)\lra\cP(\tt^*)^W.$\\
\ind The $W_K-$chamber in $\R^k=i\tt_0^*$ corresponding to chosen positive roots $\D_K^+$ is
$$
C=\{\l\in\R^k;\,\,\l_1>\l_2>\ldots>\l_{k-1}>|\l_k|>0\}.
$$
The set $\cD$ of $W-$chambers contained in $C$ consists of two elements:
$$
D_1=\{\l\in\R^k;\,\,\l_1>\l_2>\cdots>\l_{k-1}>\l_k>0\}
$$
and
$$
D_2=\{\l\in\R^k;\,\,\l_1>\l_2>\cdots>\l_{k-1}>-\l_k>0\}.
$$
The closure $\ov{D}_1$ is fundamental domain for the action of $W$ on $\R^k,$ i.e. each $W-$orbit in $\R^k$ intersects with $\ov{D}_1$ in one point. We saw that the reducibility criteria imply that $\Lambda(\s,\nu)\in\R^k$ whenever $\pi^{\s,\nu}$ is reducible. We denote by $\l(\s,\nu)$ the unique point in the intersection of $W\Lambda(\s,\nu)$ with $\ov{D}_1.$ In the following theorem without loss of generality we can suppose that $\nu\geq0,$ since $\pi^{\s,\nu}$ and $\pi^{\s,-\nu}$ have the same irreducible subquotients.
\beg{tm}\beg{\en}\item[$(i)$] $\pi^{\s,\nu}$ is reducible if and only if its infinitesimal character is $\chi_{\l}$ for some $\l\in\Lambda,$ where
$$
\beg{array}{c}
\Lambda=\le\{\l\in\Z_+^k\cup\le(\fr{1}{2}+\Z_+\ri)^k;\,\,\l_1>\l_2>\cdots>\l_{k-1}>\l_k\geq0\ri\}.
\end{array}
$$
We write $\Lambda$ as the disjoint union $\Lambda^*\cup\Lambda^0,$ where
$$
\beg{array}{c}
\Lambda^*=\le\{\l\in\Z_+^k\cup\le(\fr{1}{2}+\Z_+\ri)^k;\,\,\l_1>\l_2>\cdots>\l_{k-1}>\l_k>0\ri\},
\end{array}
$$
$$
\beg{array}{c}
\Lambda^0=\le\{\l\in\Z_+^k\cup\le(\fr{1}{2}+\Z_+\ri)^k;\,\,\l_1>\l_2>\cdots>\l_{k-1}>0,\,\,\l_k=0\ri\}.
\end{array}
$$
\item[$(ii)$] For $\l\in\Lambda^*$ there exist $k$ ordered pairs $(\s,\nu),$ $\s=(n_1,\ldots,n_{k-1})\in\hat{M},$ $\nu\geq0,$ sucha that $\chi_{\l}$ is the infinitesimal character of $\pi^{\s,\nu}.$ These ordered pairs are:
\beg{\en}
\item[$(a)$] $\nu=\l_k,$ $n_1=\l_1-k+\fr{3}{2},$ $n_2=\l_2-k+\fr{5}{2},$ $\ldots$ $n_{k-1}=\l_{k-1}-\fr{1}{2}.$
\item[$(bj)$] $\nu=\l_j,$ $n_1=\l_1-k+\fr{3}{2},$ $\ldots$ $n_{j-1}=\l_{j-1}-k+j-\fr{1}{2},$\lb $n_j=\l_{j+1}-k+j+\fr{1}{2},$ $\ldots$ $n_{k-1}=\l_k-\fr{1}{2},$ $2\leq j\leq k-1.$
\item[$(c)$] $\nu=\l_1,$ $n_1=\l_2-k+\fr{3}{2},$ $\ldots$ $n_s=\l_{s+1}-k+s+\fr{1}{2},$ $\ldots$ $n_{k-1}=\l_k-\fr{1}{2}.$
\end{\en}
\item[$(iii)$] For $\l\in\Lambda^0,$ the ordered pair $(\s,\nu),$ $\s=(n_1,\ldots,n_{k-1})\in\hat{M},$ $\nu\in\R,$ such that $\chi_{\l}$ is the infinitesimal character of $\pi^{\s,\nu},$ is unique:
$$
\beg{array}{c}
n_1=\l_1-k+\fr{3}{2},\,\,n_2=\l_2-k+\fr{5}{2},\,\ldots\,,\,\,n_{k-1}-\fr{1}{2},\qu\nu=\l_k=0.
\end{array}
$$
\end{\en}
\end{tm}

{\bf Proof:} $(i)$ We already know that for reducible elementary representation $\pi^{\s,\nu}$ one has $\Lambda(\s,\nu)\in\Z^k\cup\le(\fr{1}{2}+\Z\ri)^k.$ As the Weyl group $W=W(\gg,\tt)$ consists of all permutations of coordinates combined with multiplying some of the coordinates with $-1,$ we conclude that the infinitesimal character of $\pi^{\s,\nu}$ is $\chi_{\l}$ for some $\l\in\Lambda.$ The sufficiency will follow from the proofs of $(ii)$ and $(iii).$\\
\ind $(ii)$ Let $\l\in\Lambda^*$ and suppose that $\chi_{\l}$ is the infinitesimal character of $\pi^{\s,\nu}.$ This means that $\Lambda(\s,\nu)$ and $\l$ are $W-$conjugated. Now, since $\l_j>0,$ $\fa j\in\{1,\ldots,k\},$ $n_s-k+s+\fr{1}{2}>0,$ $\fa s\in\{1,\ldots,k-1\}$ and $\nu\geq0,$ we conclude that necessarily $\nu=\l_j$ for some $j\in\{1,\ldots,k\}.$ We inspect now each of these $k$ possibilities.\\
\ind $(a)$ $\nu=\l_k.$ Then necessarily
$$
\beg{array}{c}
n_1=\l_1-k+\fr{3}{2},\,\,\ldots\,\,,\,\,n_{k-1}=\l_{k-1}-\fr{1}{2}.
\end{array}
$$ 
We check now that so defined $(k-1)-$tuple $(n_1,\ldots,n_{k-1})$ is indeed in $\hat{M}.$ For $1\leq j\leq k-2$ we have
$$
\beg{array}{c}
n_j-n_{j+1}=\le(\l_j-k+\fr{2j+1}{2}\ri)-\le(\l_{j+1}-k+\fr{2j+3}{2}\ri)=\l_j-\l_{j+1}-1\in\Z_+.
\end{array}
$$
Further, if $\l\in\N^k,$ then $\l_{k-1}\geq2,$ thus $n_{k-1}\geq\fr{3}{2},$ and if $\l\in\le(\fr{1}{2}+\Z_+\ri)^k,$ then $\l_{k-1}\geq\fr{3}{2},$ thus $n_{k-1}\geq1.$ Especially, $\s=(n_1,\ldots,n_{k-1})\in\hat{M}.$ Finally, we see that $\nu=\l_k\leq\l_{k-1}-1=n_{k-1}-\fr{1}{2},$ so we conclude that the elementary representation $\pi^{~s,\nu}$ is reducible.\\
\ind $(bj)$ $\nu=\l_j$ for some $j\in\{2,\ldots,k-1\}.$ Then necessarily $n_1=\l_1-k+\fr{3}{2},$ $\ldots,$ $n_{j-1}=\l_{j-1}-k+\fr{2j-1}{2},n_j=\l_{j+1}-k+\fr{2j+1}{2},$ $\ldots,$ $n_{k-1}=\l_k-\fr{1}{2}.$ We check now that so defined $(k-1)-$tuple $(n_1,\ldots,n_{k-1})$ is indeed in $\hat{M}.$ For $1\leq s\leq j-2$ we have
$$
\beg{array}{c}
n_s-n_{s+1}=\le(\l_s-k+\fr{2s+1}{2}\ri)-\le(\l_{s+1}-k+\fr{2s+3}{2}\ri)=\l_s-\l_{s+1}-1\in\Z_+.
\end{array}
$$
Further,
$$
\beg{array}{c}
n_{j-1}-n_j=\le(\l_{j-1}-k+\fr{2j-1}{2}\ri)-\le(\l_{j+1}-k+\fr{2j+1}{2}\ri)=\l_{j-1}-\l_{j+1}-1\in\N.
\end{array}
$$
For $j\leq s\leq k-2$ we have
$$
\beg{array}{c}
n_s-n_{s+1}=\le(\l_{s+1}-k+\fr{2s+1}{2}\ri)-\le(\l_{s+2}-k+\fr{2s+3}{2}\ri)=\l_{s+1}-\l_{s+2}-1\in\Z_+.
\end{array}
$$
Finally, $n_{k-1}=\l_k-\fr{1}{2}\in\fr{1}{2}\Z_+.$ Thus, $\s=(n_1,\ldots,n_{k-1})\in\hat{M}.$\\
\ind We check now that the elementary representation $\pi^{\s,\nu}$ is reducible. We have
$$
\beg{array}{c}
1\leq\l_{j-1}-\l_j=n_{j-1}+k-j+\fr{1}{2}-\nu\qu\Lra\qu\nu\leq n_{j-1}+k-j-\fr{1}{2}
\end{array}
$$
and
$$
\beg{array}{c}
1\leq\l_j-\l_{j+1}=\nu-n_j-k+j+\fr{1}{2}\qu\Lra\qu\nu\geq n_j+k-j+\fr{1}{2}.
\end{array}
$$
Thus,
$$
\beg{array}{c}
\nu\in\le\{n_j+k-j+\fr{1}{2},\ldots,n_{j-1}+k-j-\fr{1}{2}\ri\}
\end{array}
$$
and we conclude that $\pi^{\s,\nu}$ is reducible.\\
\ind $(c)$ $\nu=\l_1.$ Then necessarily
$$
\beg{array}{c}
n_1=\l_2-k+\fr{3}{2},\ldots,n_s=\l_{s+1}-k+\fr{2s+1}{2},\ldots,n_{k-1}=\l_k-\fr{1}{2}.
\end{array}
$$
As before we see that for $1\leq s\leq k-2$ we have
$$
n_s-n_{s+1}=\l_{s+1}-\l_{s+2}-1\in\Z_+.
$$
Further, $n_{k-1}=\l_k-\fr{1}{2}\in\fr{1}{2}\Z_+.$ Thus, $\s=(n_1,\ldots,n_{k-1})\in\hat{M}.$ Finally,
$$
\beg{array}{c}
1\leq\l_1-\l_2=\nu-\le(n_1+k-\fr{3}{2}\ri)\qu\Lra\qu\nu\geq n_+k-\fr{1}{2},
\end{array}
$$ 
i.e.
$$
\beg{array}{c}
\nu\in\le\{n_1+k-\fr{1}{2},n_1+k+\fr{1}{2},n_1+k+\fr{3}{2},\ldots\ri\}.
\end{array}
$$
Thus, the elementary representation $\pi^{\s,\nu}$ is reducible.\\
\ind $(iii)$ Let $\l\in\Lambda^0$ and suppose that the elementary representation $\pi^{\s,\nu},$ $\s=(n_1,\ldots,n_{k-1})\in\hat{M},$ $\nu\in\R,$ has infinitesimal character $\chi_{\l}.$ As in the proof of $(ii)$ we conclude that necessarily $|\nu|=\l_j$ for some $j\in\{1,\ldots,k\}.$ The assumption $j<k$ would imply $n_{k-1}+\fr{1}{2}=\l_k=0$ and this is impossible since $n_{k-1}\geq0.$ Thus, we conclude that $j=k,$ i.e. $\nu=\l_k=0.$ It follows that
$$
\beg{array}{c}
n_1=\l_1-k+\fr{3}{2},\,\,n_2=\l_2-k+\fr{5}{2},\,\,\ldots\,\,,\,\,n_{k-1}=\l_{k-1}-\fr{1}{2}.
\end{array}
$$
As before we conclude that so defined $\s=(n_1,\ldots,n_{k-1})$ is in $\hat{M}.$ Finally, $n_j$ are nonintegral and $\nu=0,$ therefore $\pi^{\s,0}$ is reducible.\\

\ind Fix now $\l\in\Lambda^*.$ By $(ii)$ in Theorem 1. there exist $k$ pairs $(\s,\nu)\in\hat{M}\tim\fr{1}{2}\N$ such that $\chi_{\l}$ is the infinitesimal character of $\pi^{\s,\nu}.$ Denote them by $(\s_j,\nu_j),$ $1\leq j\leq k,$ where $\nu_j=\l_j$ and
$$
\beg{array}{ll}
(c)&\s_1=\le(\l_2-k+\fr{3}{2},\ldots,\l_{s+1}-k+\fr{2s+1}{2},\ldots,\l_k-\fr{1}{2}\ri),\\
&\\
(bj)&\s_j=\le(\l_1-k+\fr{3}{2},\ldots,\l_{j-1}-k+\fr{2j-1}{2},\l_{j+1}-k+\fr{2j+1}{2},\ldots,\l_k-\fr{1}{2}\ri),\\
&\qqu2\leq j\leq k-1,\\
&\\
(a)&\s_k=\le(\l_1-k+\fr{3}{2},\ldots,\l_s-k+\fr{2s+1}{2},\ldots,\l_{k-1}-\fr{1}{2}\ri).
\end{array}
$$
\ind There are altogether $k+2$ mutually infinitesimally inequivalent irreducible subquotients of the reducible elementary representations $\pi^{\s_1,\nu_1},\ldots,\pi^{\s_k,\l_k}$ which we denote by $\t_1^{\l},\ldots,\t_k^{\l},\o_+^{\l},\o_-^{\l}:$
$$
\beg{array}{l}
\t_1^{\l}=\t^{\s_1,\nu_1},\\
\\
\t_2^{\l}=\o^{\s_1,\nu_1}\cong\t^{\s_2,\nu_2},\\
\\
\vdots\\
\\
\t_j^{\l}=\o^{\s_{j-1},\nu_{j-1}}\cong\t^{\s_j,\nu_j},\\
\\
\vdots\\
\\
\t_k^{\l}=\o^{\s_{k-1},\nu_{k-1}}\cong\t^{\s_k,\nu_k},\\
\\
\o_+^{\l}=\o^{\s_k,\nu_k,+},\\
\\
\o_-^{\l}=\o^{\s_k,\nu_k,-}.
\end{array}
$$
\ind The $K-$spectra of these irreducible representations consist of all\lb $q=(m_1,\ldots,m_k)\in\hat{K}\cap\le(\l_1+\fr{1}{2}+\Z\ri)^k$ that satisfy:
$$
\beg{array}{ll}
\G(\t_1^{\l}):&\l_1-k+\fr{1}{2}\geq m_1\geq\l_2-k+\fr{3}{2}\geq\cdots\geq m_{k-1}\geq\l_k-\fr{1}{2}\geq|m_k|,\\
&\,\,\vdots\\
\G(\t_j^{\l}):&m_1\geq\l_1-k+\fr{3}{2}\geq m_2\geq\cdots\geq m_{j-1}\geq\l_{j-1}-k+j-\fr{1}{2},\\
&\l_j-k+j-\fr{1}{2}\geq m_j\geq\cdots\geq m_{k-1}\geq\l_k-\fr{1}{2}\geq|m_k|,\\
&\,\,\vdots\\
\G(\t_k^{\l}):&m_1\geq\l_1-k+\fr{3}{2}\geq m_2\geq\l_2-k+\fr{5}{2}\geq\cdots m_{k-1}\geq\l_{k-1}-\fr{1}{2},\\
&\l_k-\fr{1}{2}\geq|m_k|,\\
&\\
\G(\o_+^{\l}):&m_1\geq\l_1-k+\fr{3}{2}\geq\cdots\geq m_{k-1}\geq\l_{k-1}-\fr{1}{2}\geq m_k\geq\l_k+\fr{1}{2},\\
&\\
\G(\o_-^{\l}):&m_1\geq\l_1-k+\fr{3}{2}\geq\cdots\geq m_{k-1}\geq\l_{k-1}-\fr{1}{2}\geq-m_k\geq\l_k+\fr{1}{2}.
\end{array}
$$
\ind It is obvious that each of these representations $\pi$ has one $D_1-$corner, we denote it by $q_1(\pi),$ and one $D_2-$corner, we denote it by $q_2(\pi).$ The list is:
$$
\beg{array}{l}
q_1(\t_1^{\l})=\le(\l_2-k+\fr{3}{2},\ldots,\l_{k-1}-\fr{3}{2},\l_k-\fr{1}{2},-\l_k+\fr{1}{2}\ri),\\
q_2(\t_1^{\l})=\le(\l_2-k+\fr{3}{2},\ldots,\l_{k-1}-\fr{3}{2},\l_k-\fr{1}{2},\l_k-\fr{1}{2}\ri),\\
q_1(\t_j^{\l})=\le(\l_1-k+\fr{3}{2},\ldots,\l_{j-1}-k+j-\fr{1}{2},\l_{j+1}-k+j+\fr{1}{2},\ldots,\l_k-\fr{1}{2},-\l_k+\fr{1}{2}\ri),\\
q_2(\t_j^{\l})=\le(\l_1-k+\fr{3}{2},\ldots,\l_{j-1}-k+j-\fr{1}{2},\l_{j+1}-k+j+\fr{1}{2},\ldots,\l_k-\fr{1}{2},\l_k-\fr{1}{2}\ri),\\
q_1(\t_k^{\l})=\le(\l_1-k+\fr{3}{2},\l_2-k+\fr{5}{2},\ldots,\l_{k-1}-\fr{1}{2},-\l_k+\fr{1}{2}\ri),\\
q_2(\t_k^{\l})=\le(\l_1-k+\fr{3}{2},\l_2-k+\fr{5}{2},\ldots,\l_{k-1}-\fr{1}{2},\l_k-\fr{1}{2}\ri),\\
q_1(\o_+^{\l})=\le(\l_1-k+\fr{3}{2},\l_2-k+\fr{5}{2},\ldots,\l_{k-1}-\fr{1}{2},\l_k+\fr{1}{2}\ri),\\
q_2(\o_+^{\l})=\le(\l_1-k+\fr{3}{2},\l_2-k+\fr{5}{2},\ldots,\l_{k-1}-\fr{1}{2},\l_{k-1}-\fr{1}{2}\ri),\\
q_1(\o_-^{\l})=\le(\l_1-k+\fr{3}{2},\l_2-k+\fr{5}{2},\ldots,\l_{k-1}-\fr{1}{2},-\l_{k-1}+\fr{1}{2}\ri),\\
q_2(\o_-^{\l})=\le(\l_1-k+\fr{3}{2},\l_2-k+\fr{5}{2},\ldots,\l_{k-1}-\fr{1}{2},-\l_k-\fr{1}{2}\ri).
\end{array}
$$
We inspect now which of these corners are fundamental. Since
$$
\beg{array}{l}
\r_K-\r_P^{D_1}=\le(k-\fr{3}{2},k-\fr{5}{2},\ldots,\fr{1}{2},-\fr{1}{2}\ri),\\
\\
\r_K-\r_P^{D_2}=\le(k-\fr{3}{2},k-\fr{5}{2},\ldots,\fr{1}{2},\fr{1}{2}\ri),
\end{array}
$$
we have
$$
\beg{array}{ll}
q_1(\t_1^{\l})+\r_k-\r_P^{D_1}=(\l_2,\ldots,\l_k,-\l_k),&\mr{not\,\,fundamental},\\
q_2(\t_1^{\l})+\r_k-\r_P^{D_2}=(\l_2,\ldots,\l_k,\l_k),&\mr{not\,\,fundamental},\\
q_1(\t_j^{\l})+\r_k-\r_P^{D_1}=(\l_1,\ldots,\l_{j-1},\l_{j+1},\ldots,\l_k,-\l_k),&\mr{not\,\,fundamental},\\
q_2(\t_j^{\l})+\r_k-\r_P^{D_2}=(\l_1,\ldots,\l_{j-1},\l_{j+1},\ldots,\l_k,\l_k),&\mr{not\,\,fundamental},\\
q_1(\t_k^{\l})+\r_k-\r_P^{D_1}=(\l_1,\ldots,\l_{k-1},-\l_k),&\mr{fundamental},\\
q_2(\t_k^{\l})+\r_k-\r_P^{D_2}=(\l_1,\ldots,\l_{k-1},\l_k),&\mr{fundamental},\\
q_1(\o_+^{\l})+\r_k-\r_P^{D_1}=(\l_1,\ldots,\l_{k-1},\l_k),&\mr{fundamental},\\
q_2(\o_+^{\l})+\r_k-\r_P^{D_2}=(\l_1,\ldots,\l_{k-1},\l_{k-1}),&\mr{not\,\,fundamental},\\
q_1(\o_-^{\l})+\r_k-\r_P^{D_1}=(\l_1,\ldots,\l_{k-1},-\l_{k-1}),&\mr{not\,\,fundamental},\\
q_2(\o_-^{\l})+\r_k-\r_P^{D_2}=(\l_1,\ldots,\l_{k-1},-\l_k),&\mr{fundamental}.
\end{array}
$$
\ind Notice that finite dimensional $\t_1^{\l}$ is not unitary and $q_1(\t_1^{\l})\not=q_2(\t_1^{\l})$ unless it is the trivial $1-$dimensional representation $(\l=\le(k-\fr{1}{2},k-\fr{3}{2},\ldots,\fr{1}{2}\ri))$ when $q_1(\t_1^{\l})=q_2(\t_1^{\l})=(0,\ldots,0).$ Next, $\t_j^{\l}$ for $2\leq j\leq k$ is not unitary and $q_1(\t_j^{\l})\not=q_2(t_j^{\l}).$ Finally, $\o_+^{\l}$ and $\o_-^{\l}$ are unitary (these are the discrete series representations) and each of them has one fundamental corner; the other corner is not fundamental.\\

\ind We consider now the case $\l\in\Lambda^0,$ so $\l_k=0.$ Then the unique pair $(\s,\nu)\in\hat{M}\tim\R,$ such that $\chi_{\l}$ is the infinitesimal character of the elementary representation $\pi^{\s,\nu},$ is
$$
\beg{array}{ll}
\s=\le(\l_1-k+\fr{3}{2},\l_2-k+\fr{5}{2},\ldots,\l_{k-1}-\fr{1}{2}\ri),&\nu=0.
\end{array}
$$
The elementary representation $\pi^{\s,0}$ is unitary and it is direct sum of two unitary irreducible representations $\o_+^{\l}$ and $\o_-^{\l}.$ Their $K-$spectra consist of all $q=(m_1,\ldots,m_k)\in\hat{K}\cap\le(\fr{1}{2}+\Z\ri)^k$ that satisfy
$$
\beg{array}{ll}
\G(\o_+^{\l}):&m_1\geq\l_1-k+\fr{3}{2}\geq\cdots\geq m_{k-1}\geq\l_{k-1}-\fr{1}{2}\geq m_k\geq\fr{1}{2},\\
&\\
\G(\o_-^{\l}):&m_1\geq\l_1-k+\fr{3}{2}\geq\cdots\geq m_{k-1}\geq\l_{k-1}-\fr{1}{2}\geq-m_k\geq\fr{1}{2}.
\end{array}
$$
Again each of these representations have one $D_1-$corner and one $D_2-$corner:
$$
\beg{array}{l}
q_1(\o_+^{\l})=\le(\l_1-k+\fr{3}{2},\l_2-k+\fr{5}{2},\ldots,\l_{k-1}-\fr{1}{2},\fr{1}{2}\ri),\\
q_2(\o_+^{\l})=\le(\l_1-k+\fr{3}{2},\l_2-k+\fr{5}{2},\ldots,\l_{k-1}-\fr{1}{2},\l_{k-1}-\fr{1}{2}\ri),\\
q_1(\o_-^{\l})=\le(\l_1-k+\fr{3}{2},\l_2-k+\fr{5}{2},\ldots,\l_{k-1}-\fr{1}{2},-\l_{k-1}+\fr{1}{2}\ri),\\
q_2(\o_-^{\l})=\le(\l_1-k+\fr{3}{2},\l_2-k+\fr{5}{2},\ldots,\l_{k-1}-\fr{1}{2},-\fr{1}{2}\ri).
\end{array}
$$
Two of them are fundamental:
$$
\beg{array}{ll}
q_1(\o_+^{\l})+\r_k-\r_P^{D_1}=(\l_1,\ldots,\l_{k-1},0),&\mr{fundamental},\\
q_2(\o_+^{\l})+\r_k-\r_P^{D_2}=(\l_1,\ldots,\l_{k-1},\l_{k-1}),&\mr{not\,\,fundamental},\\
q_1(\o_-^{\l})+\r_k-\r_P^{D_1}=(\l_1,\ldots,\l_{k-1},-\l_{k-1}),&\mr{not\,\,fundamental},\\
q_2(\o_-^{\l})+\r_k-\r_P^{D_2}=(\l_1,\ldots,\l_{k-1},0),&\mr{fundamental}.
\end{array}
$$
Thus, we see that again each of these unitary representation has one fundamental corner and the other corner is not fundamental.\\

\ind To sumarize, we see that $\pi\in\wp{G}^0$ with exactly one fundamental corner is unitary; its fundamental corner we denote by $q(\pi).$ For all the others $\pi\in\hat{G}^0$ one has $q_1(\pi)=q_2(\pi)$ and we denote by $q(\pi)$ this unique corner of $\pi.$
\beg{tm} $\pi\map q(\pi)$ is a bijection of $\hat{G}^0$ onto $\hat{K}.$
\end{tm}

{\bf Proof:} We have
$$
\hat{G}^0=\bigcup_{j=1}^{k}\{\t_j^{\l};\,\,\l\in\Lambda_j^*\}\cup\{\o_+^{\l};\,\,\l\in\Lambda\}\cup\{\o_-^{\l};\,\,\l\in\Lambda\},
$$
where $\Lambda_1^*=\le\{\le(k-\fr{1}{2},k-\fr{3}{2},\ldots,\fr{1}{2}\ri)\ri\}$ and for $2\leq j\leq k$
$$
\beg{array}{c}
\Lambda_j^*=\le\{\l\in\Lambda^*;\,\,\l_{j-1}>k-j+\fr{3}{2}\,\,\mr{and}\,\,\l_s=k-s+\fr{1}{2}\,\,\mr{for}\,\,j\leq s\leq k\ri\}.
\end{array}
$$
\ind Let $q=(m_1,\ldots,m_k)\in\hat{K}.$ We have three possibilities:\\
\ind$(1)$ $m_k=0.$ Then $q\in\Z_+^k$ and $m_1\geq m_2\geq\cdots\geq m_{k-1}\geq0.$ Let
$$
j=\min\,\{s;\,\,1\leq s\leq k,\,\,m_s=0\}.
$$
Set
$$
\beg{array}{ll}
\l_s=m_s+k-s+\fr{1}{2},&\qu1\leq s\leq j-1,\\
&\\
\l_s=k-s+\fr{1}{2},&\qu j\leq s\leq k.
\end{array}
$$
Then for $1\leq s\leq j-2$ we have $\l_s-\l_{s+1}=m_s-m_{s+1}+1\geq1,$ next $\l_{j-1}-\l_j=m_{j-1}+1\geq2,$ further, for $j\leq s\leq k-1$ we have $\l_s-\l_{s+1}=1,$ and finally $\l_k=\fr{1}{2}.$ Thus, we see that $\l\in\Lambda^*.$ If $j=1$ we see that $\l$ is the unique element of $\Lambda_1^*.$ If $j\geq2$ we have $m_{j-1}>0$ and so
$$
\beg{array}{c}
\l_{j-1}=m_{j-1}+k-j+1+\fr{1}{2}>k-j+\fr{3}{2}
\end{array}
$$
i.e. $\l\in\Lambda_j^*.$ From the definition of $\l$ we see that $q=q(\t_j^{\l}).$\\
\ind$(2)$ $m_k>0.$ Set now
$$
\beg{array}{c}
\l_j=m_j+k-j-\fr{1}{2}.
\end{array}
$$
Then $\l_j-\l_{j+1}=m_j-m_{j-1}+1\geq1$ for $1\leq j\leq k-1$ and $\l_k=m_k-\fr{1}{2}\geq0.$ Thus, $\l\in\Lambda$ and one sees that $q=q(\o_+^{\l}).$\\
\ind$(3)$ $m_k<0.$ Set now
$$
\beg{array}{c}
\l_j=m_j+k-j-\fr{1}{2},\,\,1\leq j\leq k-1,\qu\l_k=-m_k-\fr{1}{2}.
\end{array}
$$
Then $\l_j-\l_{j+1}=m_j-m_{j+1}+1\geq1$ for $1\leq j\leq k-2,$ further\lb $\l_{k-1}-\l_k=m_{k-1}+m_k+1=m_{k-1}-|m_k|+1\geq1,$ and finally $\l_k=|m_k|-\fr{1}{2}\geq0.$ Thus, $\l\in\Lambda$ and one sees that $q=q(\o_-^{\l}).$\\
\ind We have proved that $\pi\map q(\pi)$ is a surjection of $\hat{G}^0$ onto $\hat{K}.$ From the proof we see that this map is injective too.\\

\ind Consider now minimal $K-$types in the sense of Vogan: we say that $q\in\hat{K}$ is a {\bf minimal} $K-${\bf type} of the representation $\pi$ if $q\in\G(\pi)$ and
$$
\|q+2\r_K\|=\min\,\{\|q^{\pr}+2\r_K\|;\,\,q^{\pr}\in\G(\pi)\}.
$$
For $q\in\hat{K}$ we have
$$
\|q+2\r_K\|^2=(m_1+2k-2)^2+(m_2+2k-4)^2+\cdots+(m_{k-1}+2)^2+m_k^2
$$
and so we find:\\
\ind If $\l\in\Lambda\cap\le(\fr{1}{2}+\Z\ri)^k,$ i.e. $\l\in\Lambda^*$ and $\G(\t_j^{\l})\sub\Z^k,$ the representation $\t_j^{\l}$ has one minimal $K-$type which we denote by $q^V(\t_j^{\l}):$
$$
\beg{array}{l}
q^V(\t_1^{\l})=\le(\l_2-k+\fr{3}{2},\l_3-k+\fr{5}{2},\ldots,\l_k-\fr{1}{2},0\ri),\\
\\
q^V(\t_j^{\l})=\le(\l_1-k+\fr{3}{2},\ldots,\l_{j-1}-k+j-\fr{1}{2},\l_{j+1}-k+j+\fr{1}{2},\ldots,\l_k-\fr{1}{2},0\ri),\\
\qqu\qqu2\leq j\leq k-1,\\
\\
q^V(\t_k^{\l})=\le(\l_1-k+\fr{3}{2},\l_2-k+\fr{5}{2},\ldots,\l_{k-1}-\fr{1}{2},0\ri).
\end{array}
$$

\newpage

\ind If $\l\in\Lambda\cap\Z^k,$ i.e. $\G(\t_j^{\l})\sub\le(\fr{1}{2}+\Z\ri)^k,$ the representation $\t_j^{\l}$ has two minimal $K-$types $q_1^V(\t_j^{\l})$ and $q_2^V(\t_j^{\l}):$
$$
\beg{array}{l}
q_1^V(\t_1^{\l})=\le(\l_2-k+\fr{3}{2},\l_3-k+\fr{5}{2},\ldots,\l_k-\fr{1}{2},\fr{1}{2}\ri),\\
\\
q_2^V(\t_1^{\l})=\le(\l_2-k+\fr{3}{2},\l_3-k+\fr{5}{2},\ldots,\l_k-\fr{1}{2},-\fr{1}{2}\ri),\\
\\
q_1^V(\t_j^{\l})=\le(\l_1-k+\fr{3}{2},\ldots,\l_{j-1}-k+j-\fr{1}{2},\l_{j+1}-k+j+\fr{1}{2},\ldots,\l_k-\fr{1}{2},\fr{1}{2}\ri),\\
\qqu\qqu2\leq j\leq k-1,\\
\\
q_2^V(\t_j^{\l})=\le(\l_1-k+\fr{3}{2},\ldots,\l_{j-1}-k+j-\fr{1}{2},\l_{j+1}-k+j+\fr{1}{2},\ldots,\l_k-\fr{1}{2},-\fr{1}{2}\ri),\\
\qqu\qqu2\leq j\leq k-1,\\
\\
q_1^V(\t_k^{\l})=\le(\l_1-k+\fr{3}{2},\l_2-k+\fr{5}{2},\ldots,\l_{k-1}-\fr{1}{2},\fr{1}{2}\ri),\\
\\
q_2^V(\t_k^{\l})=\le(\l_1-k+\fr{3}{2},\l_2-k+\fr{5}{2},\ldots,\l_{k-1}-\fr{1}{2},-\fr{1}{2}\ri).
\end{array}
$$
\ind Finally, for every $\l\in\Lambda$ the representation $\o_{\pm}^{\l}$ has one minimal $K-$type $q^V(\o_{\pm}^{\l}):$
$$
\beg{array}{l}
q^V(\o_+^{\l})=\le(\l_1-k+\fr{3}{2},\l_2-k+\fr{5}{2},\ldots,\l_{k-1}-\fr{1}{2},\l_k+\fr{1}{2}\ri),\\
\\
q^V(\o_-^{\l})=\le(\l_1-k+\fr{3}{2},\l_2-k+\fr{5}{2},\ldots,\l_{k-1}-\fr{1}{2},-\l_k-\fr{1}{2}\ri).
\end{array}
$$
\ind So we see that if $\pi\in\wp{G}^0$ has two minimal $K-$types it is not unitary. Further, every $\pi\in\hat{G}^0$ has one minimal $K-$type $q^V(\pi)$ and it coincides with $q(\pi).$ But there exist nonunitary representations in $\wp{G}^0$ that have one mi\-ni\-mal $K-$type: this property have all $\t_j^{\l}$ for $\l\in\Lambda\cap\le(\fr{1}{2}+\Z_+\ri)^k$ that are not subquotinets of the ends of compelemntary series. In other words, unitarity of a representation $\pi\in\wp{G}^0$ is not characterized by having unique minimal $K-$type.

\section{Representations of $\Spin(2k+1,1)$}

\ind Now $M=\Spin(2k).$ Cartan subalgebra $\tt_0$ of $\kk_0$ (resp. $\tt$ of $\kk)$ is also Cartan subalgebra of $\mm_0$ (resp. $\mm).$ The root systems are:
$$
\D_K=\D(\kk,\tt)=\{\pm\a_p\pm\a_q;\,\,1\leq p,q\leq k,\,\,p\not=q\}\cup\{\pm\a_p;\,\,1\leq p\leq k\}
$$
and
$$
\D_M=\D(\mm,\tt)=\{\pm\a_p\pm\a_q;\,\,1\leq p,q\leq k,\,\,p\not=q\}.
$$
We choose positive roots:
$$
\D_K^+=\{\a_p\pm\a_q;\,\,1\leq p<q\leq k\}\cup\{\a_p;\,\,1\leq p\leq k\},
$$
$$
\D_M^+=\{\a_p\pm\a_q;\,\,1\leq p<q\leq k\}.
$$
The corresponding Weyl chambers in $\R^k=i\tt_0^*$ are
$$
C_K=\{\l\in\R^k;\,\,\l_1>\l_2>\cdots>\l_{k-1}>\l_k>0\}
$$
with the closure
$$
\ov{C}_K=\{\l\in\R^k;\,\,\l_1\geq\l_2\geq\cdots\geq\l_{k-1}\geq\l_k\geq0\}
$$
and
$$
C_M=\{\l\in\R^k;\,\,\l_1>\l_2>\cdots>\l_{k-1}>|\l_k|>0\}
$$
with the closure
$$
\ov{C}_M=\{\l\in\R^k;\,\,\l_1\geq\l_2\geq\cdots\geq\l_{k-1}\geq|\l_k|\}.
$$
The halfsums of positive roots are
$$
\beg{array}{ccc}
\r_K=\le(k-\fr{1}{2},k-\fr{3}{2},\ldots,\fr{3}{2},\fr{1}{2}\ri)&\mr{and}&\d_{\mm}=(k-1,k-2,\ldots,1,0).
\end{array}
$$
Now
$$
\beg{array}{c}
\hat{K}=\le\{(m_1,\ldots,m_k)\in\Z_+^k\cup\le(\fr{1}{2}+\Z_+\ri)^k;\,\,m_1\geq m_2\geq\cdots\geq m_{k-1}\geq m_k\geq0\ri\}\\
\\
\hat{M}=\le\{(n_1,\ldots,n_k)\in\Z^k\cup\le(\fr{1}{2}+\Z\ri)^k;\,\,n_1\geq n_2\geq\cdots\geq n_{k-1}\geq|n_k|\ri\}.
\end{array}
$$
The branching rule is
$$
(m_1,\ldots,m_k)|M=\bigoplus_{(n_1,\ldots,n_k)\prec(m_1,\ldots,m_k)}(n_1,\ldots,n_k)
$$
where $(n_1,\ldots,n_k)\prec(m_1,\ldots,m_k)$ means that $(m_1,\ldots,m_k)\in(n_1+\Z)^k$ and
$$
m_1\geq n_1\geq m_2\geq n_2\geq\cdots\geq m_{k-1}\geq n_{k-1}\geq m_k\geq|n_k|.
$$
So by the Frobenius Reciprocity Theorem for $\s=(n_1,\ldots,n_k)\in\hat{M}$ and $\nu\in\C=\aa^*$ we have
$$
\pi^{\s,\nu}|K=\bigoplus_{(n_1,\ldots,n_k)\prec(m_1,\ldots,m_k)}(m_1,\ldots,m_k).
$$
We identify the dual $\hh^*$ with $\C^{k+1}$ so that $\l\in\hh^*$ is identified with the\lb $(k+1)-$tuple $(\l(H_1),\ldots,\l(H_k),\l(H))$ and $\tt^*=\C^k$ is identified with the subspace of $\hh^*=\C^{k+1}$ of all $(k+1)-$tuples with $0$ at the end. The infinitesimal character of the elementary representation $\pi^{\s,\nu}$ is equal $\chi_{\Lambda(\s,\nu)},$ where $\Lambda(\s,\nu)\in\hh^*$ is defined by
$$
\Lambda(\s,\nu)|\tt=\l_{\s}+\d_{\mm}\qu\mr{and}\qu\Lambda(\s,\nu)|\aa=\nu.
$$
Here $\l_{\s}$ is the highest weight of $\s$ with respect to $\D_M^+.$ Thus
$$
\Lambda(\s,\nu)=(n_1+k-1,n_2+k-2,\ldots,n_{k-1}+1,n_k,\nu).
$$
For $\s=(n_1,\ldots,n_k)\in\hat{M}\cap\Z^k$ and $\nu\in\C$ the elementary representation $\pi^{\s,\nu}$ is ireducible if and only if either $\nu\not\in\Z$ or
$$
\nu\in\{0,\pm1,\ldots,\pm|n_k|,\pm(n_{k-1}+1),\pm(n_{k-2}+2),\ldots,\pm(n_1+k-1)\}.
$$
For $\s\in\hat{M}\cap\le(\fr{1}{2}+\Z\ri)^k$ and $\nu\in\C$ the representation $\pi^{\s,\nu}$ is irreducible if and only if either $\nu\not\in\le(\fr{1}{2}+\Z\ri)$ or
$$
\beg{array}{c}
\nu\in\le\{\pm\fr{1}{2},\ldots,\pm|n_k|,\pm(n_{k-1}+1),\pm(n_{k-2}+2),\ldots,\pm(n_1+k-1)\ri\}.
\end{array}
$$
If the elementary representation $\pi^{\s,\nu}$ is reducible, it always has two irreducible subquotients which will be denoted by $\t^{\s,\nu}$ and $\o^{\s,\nu}.$ The $K-$spectra of these representations consist of all $q=(m_1,\ldots,m_k)\in\hat{K}\cap(n_1+\Z)^k$ that satisfy:
$$
\beg{array}{l}
\bullet\,\mr{If}\,\,n_{k-1}>|n_k|\,\,\mr{and}\,\,\nu\in\{\pm(|n_k|+1),\pm(|n_k|+2),\ldots,\pm n_{k-1}\}:\\
\G(\t^{\s,\nu}):\,\,m_1\geq n_1\geq\cdots\geq m_{k-1}\geq n_{k-1}\,\,\mr{and}\,\,|\nu|-1\geq m_k\geq|n_k|,\\
\G(\o^{\s,\nu}):\,\,m_1\geq n_1\geq\cdots\geq m_{k-1}\geq n_{k-1}\geq m_k\geq|\nu|.\\
\bullet\,\mr{If}\,\,n_{j-1}>n_j\,\,\mr{for\,\,some}\,\,j\in\{2,\ldots,k-1\}\,\,\mr{and}\\
\nu\in\{\pm(n_j+k-j+1),\pm(n_j+k-j+2),\ldots,\pm(n_{j-1}+k-j)\}:\\
\G(\t^{\s,\nu}):\,\,m_1\geq n_1\geq\cdots\geq m_{j-1}\geq n_{j-1}\,\,\mr{and}\\
\qqu\qqu|\nu|-k+j-1\geq m_j\geq n_j\geq\cdots\geq m_k\geq|n_k|,\\
\G(\o^{\s,\nu}):\,\,m_1\geq n_1\geq\cdots\geq m_{j-1}\geq n_{j-1}\geq m_j\geq|\nu|-k+j\,\,\mr{and}\\
\qqu\qqu n_j\geq m_{j+1}\geq\cdots\geq m_k\geq|n_k|.\\
\bullet\,\mr{If}\,\,\nu\in\{\pm(n_1+k),\pm(n_1+k+1),\ldots\}:\\
\G(\t^{\s,\nu}):\,\,|\nu|-k\geq m_1\geq n_1\geq\cdots\geq m_k\geq|n_k|,\\
\G(\o^{\s,\nu}):\,\,m_1\geq|\nu|-k+1\,\,\mr{and}\,\,n_1\geq m_2\geq n_2\geq\cdots\geq m_k\geq|n_k|.
\end{array}
$$
\ind Similarly to the preceeding case of even $n=2k$ we now write down the infinitesimal characters of reducible elementary representations $\pi^{\s,\nu}$ (and so of its irreducible subquotients $\t^{\s,\nu}$ and $\o^{\s,\nu}$ too). We know that the infinitesimal character of $\pi^{\s,\nu}$ is $\chi_{\Lambda(\s,\nu)},$ where
$$
\Lambda(\s,\nu)=(n_1+k-1,n_2+k-2,\ldots,n_{k-1}+1,n_k,\nu).
$$
Since $\nu\in\fr{1}{2}\Z\subset\R=\aa_0^*$ we have $\Lambda(\s,\nu)\in i\tt_0^*\oplus\aa_0^*=\R^{k+1}.$\\
\ind The root system of the pair $(\gg,\hh)$ is
$$
\D=\{\pm\a_p\pm\a_q;\,\,1\leq p,q\leq k+1,\,\,p\not=q\}.
$$
We choose positive roots:
$$
\D^+=\{\a_p\pm\a_q;\,\,1\leq p<q\leq k+1\}.
$$
The corresponding Weyl chamber in $\R^{k+1}$ is
$$
D=\{\l\in\R^{k+1};\,\,\l_1>\l_2>\cdots>\l_k>|\l_{k+1}|>0\}
$$
with the closure
$$
\ov{D}=\{\l\in\R^{k+1};\,\,\l_1\geq\l_2\geq\cdots\geq\l_k\geq|\l_{k+1}|\}.
$$
The Weyl group $W$ of $\D$ consists of all permutations of coordinates in $\C^{k+1}=\hh^*$ combined with multiplying even number of coordinates with $-1.$ By Harish$-$Chandra's theorem $\chi_{\l}=\chi_{\l^{\pr}}$ if and only if $\l,\l^{\pr}\in\hh^*$ are in the same $W-$orbit. As $\ov{D}$ is a fundamental domain for the action of $W$ on $\R^{k+1}=i\tt_0^*\oplus\aa_0^*,$ there exists unique $\l(\s,\nu)\in\ov{D}$ such that $\chi_{\Lambda(\s,\nu)}=\chi_{\l(\s,\nu)}.$ We now write down $\l(\s,\nu)$ for all reducible elementary representations $\pi^{\s,\nu}.$ In the following for $\s=(n_1,\ldots,n_k)\in\hat{M}$ we write $-\s$ for its contragredient class in $\hat{M}:$ $-\s=(n_1,\ldots,n_{k-1},-n_k).$ Without loss of geq-ne\-ra\-li\-ty we can suppose that $\nu\geq0$ because $\pi^{\s,\nu}$ and $\pi^{-\s,-\nu}$ have equivalent irreducible subquotients and because $\Lambda(\s,\nu)$ is $W-$conjugated with $\Lambda(-\s,-\nu):$ multiplying the last two coordinates by $-1.$\\
\ind $\bullet$ If $n_{k-1}>|n_k|$ and $\nu\in\{|n_k|+1,|n_k|+2,\ldots,n_{k-1}\}$ we have $n_{k-1}>\nu>|n_k|$ and so
$$
\l(\s,\nu)=(n_1+k-1,n_2+k-2,\ldots,n_{k-1}+1,\nu,n_k).
$$
\ind $\bullet$ If $2\leq j\leq k-1,$ $n_{j-1}>n_j$ and $\nu\in\{n_j+k-j+1,\ldots,n_{j-1}+k-j\}$ we have $n_{j-1}+k-j+1>\nu>n_j+k-j$ and so
$$
\l(\s,\nu)=(n_1+k-1,\ldots,n_{j-1}+k-j+1,\,\nu\,,n_j+k-j,\ldots,n_{k-1}+1,n_k).
$$
\ind $\bullet$ If $\nu\in\{n_1+k,n_1+k+1,\ldots\}$ we have $\nu>n_1+k-1$ and so
$$
\l(\s,\nu)=(\nu,n_1+k-1,\ldots,n_{k-1}+1,n_k).
$$
\ind Similarly to the preceeding case of even $n=2k,$ we see that now every reducible elementary representation has infinitesimal character $\chi_{\l}$ with $\l\in\Lambda,$\lb where
$$
\beg{array}{c}
\Lambda=\le\{\l\in\Z^{k+1}\cup\le(\fr{1}{2}+\Z\ri)^{k+1};\,\,\l_1>\l_2>\cdots>\l_k>|\l_{k+1}|\ri\}.
\end{array}
$$
We again write $\Lambda$ as the disjoint union $\Lambda=\Lambda^*\cup\Lambda^0,$ where
$$
\beg{array}{c}
\Lambda^*=\le\{\l\in\Z^{k+1}\cup\le(\fr{1}{2}+\Z\ri)^{k+1};\,\,\l_1\l_2>\cdots>\l_k>|\l_{k+1}|>0\ri\},\\
\\
\Lambda^0=\{\l\in\Z_+^{k+1};\,\,\l_1>\l_2>\cdots>\l_k>0,\,\,\l_{k+1}=0\}.
\end{array}
$$
\beg{tm} $(i)$ For every $\l\in\Lambda^*$ there exist $k+1$ ordered pairs $(\s,\nu),$ $\s=(n_1,\ldots,n_k)\in\hat{M},$ $\nu\geq0,$ such that $\chi_{\l}$ is the infinitesimal character of $\pi^{\s,\nu}.$ These are $(\s_j,\nu_j),$ where $\nu_j=\l_j$ for $1\leq j\leq k,$ $\nu_{k+1}=|\l_{k+1}|$ and
$$
\beg{array}{l}
\s_1=(\l_2-k+1,\l_3-k+2,\ldots,\l_k-1,\l_{k+1}),\\
\s_j=(\l_1-k+1,\ldots,\l_{j-1}-k+j-1,\l_{j+1}-k+j,\ldots,\l_k-1,\l_{k+1}),\\
\qqu2\leq j\leq k-1,\\
\s_k=(\l_1-k+1,\l_2-k+2,\ldots,\l_{k-1}-1,\l_{k+1}),\\
\s_{k+1}=\le\{\beg{array}{ll}
\!\!(\l_1-k+1,\l_2-k+2,\ldots,\l_{k-1}-1,\l_k)&\,\,\mr{if}\,\,\l_{k+1}>0\\
\!\!(\l_1-k+1,\l_2-k+2,\ldots,\l_{k-1}-1,-\l_k)&\,\,\mr{if}\,\,\l_{k+1}<0
\end{array}
\ri.
\end{array}
$$
$\pi^{\s_j,\nu_j},$ $1\leq j\leq k,$ are reducible, while $\pi^{\s_{k+1},\nu_{k+1}}$ is irreducible.\\
\ind $(ii)$ For $\l\in\Lambda^0$ there exist $k+2$ ordered pairs $(\s,\nu),$ $\s=(n_1,\ldots,n_k)\in\hat{M},$ $\nu\geq0,$ such that $\chi_{\l}$ is the infinitesimal character of $\pi^{\s,\nu}.$ These are the $(\s_j,\nu_j),$ where $\nu_j=\l_j$ for $1\leq j\leq k,$ $\nu_{k+1}=\nu_{k+2}=0$ and
$$
\beg{array}{l}
\s_1=(\l_2-k+1,\l_3-k+2,\ldots,\l_k-1,0),\\
\s_j=(\l_1-k+1,\ldots,\l_{j-1}-k+j-1,\l_{j+1}-k+j,\ldots,\l_k-1,0),\\
\qqu2\leq j\leq k-1,\\
\s_k=(\l_1-k+1,\l_2-k+2,\ldots,\l_{k-1}-1,0),\\
\s_{k+1}=(\l_1-k+1,\l_2-k+2,\ldots,\l_{k-1}-1,\l_k),\\
\s_{k+2}=(\l_1-k+1,\l_2-k+2,\ldots,\l_{k-1}-1,-\l_k).
\end{array}
$$
$\pi^{\s_j,\nu_j},$ $1\leq j\leq k,$ are reducible, while $\pi^{\s_{k+1},\nu_{k+1}}$ and $\pi^{\s_{k+2},\nu_{k+2}}$ are irreducible.
\end{tm}

\newpage

{\bf Proof:} $(i)$ Fix $\l\in\Lambda^*$ and let $(\s,\nu),$ $\s=(n_1,\ldots,n_k)\in\hat{M},$ $\nu\geq0,$ be such that $\chi_{\l}$ is the infinitesimal character of $\pi^{\s,\nu}.$ Then $\Lambda(\s,\nu)$ and $\l$ are in the same $W-$orbit. Since $\nu\geq0$ we have necessarily $\nu=\l_j$ for some $j\leq k$ or $\nu=|\l_{k+1}|.$\\
\ind Suppose $\nu=\l_j$ for some $j\leq k.$ Since $W$ acts as permutations of coordinates combined with multiplying evene number of coordinates by $-1,$ the inequalities
$$
n_1+k-1>n_2+k-2>\cdots>n_{k-1}+1>|n_k|
$$
and
$$
\l_1>\l_2>\cdots>\l_{j-1}>\l_{j+1}>\cdots>\l_k>|\l_{k+1}|>0
$$
imply
$$
n_1+k-1=\l_1,\ldots,n_{j-1}+k-j+1=\l_{j-1},
$$
$$
n_j+k-j=\l_{j+1},\ldots,n_{k-1}+1=\l_k,\,\,n_k=\l_{k+1}.
$$
The following possibilities follow:
$$
\beg{array}{ll}
\nu_1=\l_1,&\!\!\s_1=(\l_2-k+1,\l_3-k+2,\ldots,\l_k-1,\l_{k+1}),\\
&\\
\nu_j=\l_j,&\!\!\s_j=(\l_1-k+1,\ldots,\l_{j-1}-k+j-1,\l_{j+1}-k+j,\ldots,\l_k-1,\l_{k+1}),\\
&\qqu2\leq j\leq k-1,\\
&\\
\nu_k=\l_k,&\!\!\s_k=(\l_1-k+1,\ldots,\l_{k-1}-1,\l_{k+1}).
\end{array}
$$
One easily checks that so defined $\s_1,\ldots,\s_k$ are really in $\hat{M}$ and that $\pi^{\s_j,\nu_j}$ are reducible.\\
\ind Suppose now that $\l_{k+1}>0$ and $\nu=\l_{k+1}.$ Then it follows that necessarily
$$
n_1+k-1=\l_1,n_2+k-2=\l_2,\ldots,n_{k-1}+1=\l_{k-1},n_k=\l_k,
$$
i.e.
$$
n_1=\l_1-k+1,n_2=\l_2-k+2,\ldots,n_{k-1}=\l_{k-1}-1,n_k=\l_k.
$$
On the other hand, if $\l_{k+1}<0,$ hence $\nu=-\l_{k+1},$ we see that in $W-$action which $\Lambda(\s,\nu)$ transforms into $\l$ there should be one more change of sign and so necessarily $n_k=-\l_k.$ Thus we have
$$
n_1=\l_1-k+1,n_2=\l_2-k+2,\ldots,n_{k-1}=\l_{k-1}-1,n_k=-\l_k.
$$
One checks that so defined
$$
\s_{k+1}=(n_1,\ldots,n_k)=(\l_1-k+1,\ldots,\l_{k-1}-1,\pm\l_k)
$$
is really in $\hat{M}.$ Further, we have $|n_k|-\nu_{k+1}=\l_k-|\l_{k+1}|\in\N.$ Thus, either $\nu_{k+1}\in\{0,1,\ldots,|n_k|-1\}$ or $\nu_{k+1}\in\{\fr{1}{2},\fr{3}{2},\ldots,|n_k|-1\}.$ Therefore, the elementary representation $\pi^{\s_{k+1},\nu_{k+1}}$ is irreducible.\\
\ind $(ii)$ Let $\l\in\Lambda^0$ and let $(\s,\nu),$ $\s=(n_1,\ldots,n_k)\in\hat{M},$ $\nu\geq0,$ be such that $\chi_{\l}$ is the infinitesimal character of $\pi^{\s,\nu}.$ As in the proof of $(i)$ we find that necessarily $\nu=\l_j$ for some $j.$ The rest of the proof for $j\leq k$ is completely the same as in $(i).$ So we are left with the case $\nu=\l_{k+1}=0.$ As in $(i)$ besause of the inequalities $n_1+k-1>n_2+k-2>\cdots>n_{k-1}+1>|n_k|$ and $\l_1>\l_2>\cdots>\l_k>0$ we get two possibilies for $\s,$ $\s=\s_{k+1}$ and $\s_{k+2}$ from the statement $(ii).$ Finally, as in the proof of $(i)$ we check that $\s_{k+1},\s_{k+2}\in\hat{M}$ and that the representations $\pi^{\s_{k+1},0}$ and $\pi^{\s_{k+2},0}$ are irreducible.\\

\ind We note that in fact the representations $\pi^{\s_{k+1},0}$ and $\pi^{\s_{k+2},0}$ are equivalent, but this is unimportant for studying and parametrizing $\wp{G}^0$ and $\hat{G}^0.$\\

\ind Fix $\l\in\Lambda.$ By Theorem 3. there exist $k$ ordered pairs $(\s,\nu),$ $\s\in\hat{M},$ $\nu\geq0,$ with reducible $\pi^{\s,\nu}$ having $\chi_{\l}$ as the infinitesimal character. There are $k+1$ mutually inequivalent irreducible subquotients of these elementary representations; we denote them $\t_1^{\l},\ldots,\t_k^{\l},\o^{\l}:$
$$
\beg{array}{l}
\t_1^{\l}=\t^{\s_1,\nu_1},\\
\\
\t_2^{\l}=\o^{\s_1,\nu_1}\cong\t^{\s_2,\nu_2},\\
\\
\vdots\\
\\
\t_j^{\l}=\o^{\s_{j-1},\nu_{j-1}}\cong\t^{\s_j,\nu_j},\\
\\
\vdots\\
\\
\t_k^{\l}=\o^{\s_{k-1},\nu_{k-1}}\cong\t^{\s_k,\nu_k},\\
\\
\o^{\l}=\o^{\s_k,\nu_k}.
\end{array}
$$
Their $K-$spectra consist of all $q=(m_1,\ldots,m_k)\in\hat{K}\cap(n_1+\Z)^k$ satisfying:
$$
\beg{array}{ll}
\G(\t_1^{\l}):&\l_1-k\geq m_1\geq\l_2-k+1\geq m_2\geq\cdots\geq\l_k-1\geq m_k\geq|\l_{k+1}|.\\
&\\
\G(\t_j^{\l}):&m_1\geq\l_1-k+1\geq\cdots\geq m_{j-1}\geq\l_{j-1}-k+j-1\,\,\mr{and}\\
&\l_j-k+j-1\geq m_j\geq\cdots\geq\l_k-1\geq m_k\geq|\l_{k+1}\,\,\mr{for}\,\,2\leq j\leq k.\\
&\\
\G(\o^{\l}):&m_1\geq\l_1-k+1\geq\cdots\geq m_{k-1}\geq\l_{k-1}-1\geq m_k\geq\l_k.
\end{array}
$$

\ind The definitons of corners and fundamental corners do not have sense when $\rank\,\kk<\rank\,\gg.$ Consider the Vogan's minimal $K-$types. Note that
$$
\|q+2\r_K\|^2=(m_1+2k-1)^2+(m_2+2k-3^2+\cdots+(m_k+1)^2,
$$
so every $\pi\in\wp{G}^0$ has unique minimal $K-$type that will be denoted by $q^V(\pi):$ this is the element $(m_1,\ldots,m_k)\in\G(\pi)$ whose every coordinate $m_j$ is the smallest possible.
\beg{tm} The map $\pi\map q^V(\pi)$ is a surjection of $\wp{G}^0$ onto $\hat{K}.$ More precisely, for $q=(m_1,\ldots,m_k)\in\hat{K}:$\\
\ind $(a)$ There exist infinitely many $\l$'s in $\Lambda$ such that $q^V(\t_1^{\l})=q.$\\
\ind $(b)$ Let $j\in\{2,\ldots,k\}.$ The number of mutually different $\l$'s $\Lambda$ such the $q^V(\t_j^{\l})=q$ is equal:
$$
\beg{array}{cl}
0&\qu\mr{if}\,\,\,m_{j-1}=m_j,\\
m_{j-1}-m_j&\qu\mr{if}\,\,\,m_{j-1}>m_j\,\,\,\mr{and}\,\,\,m_k=0,\\
2(m_{j-1}-m_j)&\qu\mr{if}\,\,\,m_{j-1}>m_j\,\,\,\mr{and}\,\,\,m_k>0.
\end{array}
$$
\ind $(c)$ The number of $\l$'s in $\Lambda$ such that $q^V(\o^{\l})=q$ is equal:
$$
\beg{array}{cl}
0&\qu\mr{if}\,\,\,m_k=0\,\,\,\mr{or}\,\,\,m_k=\fr{1}{2},\\
1&\qu\mr{if}\,\,\,m_k=1,\\
2\le[m_k-\fr{1}{2}\ri]&\qu\mr{if}\,\,\,m_k>1.
\end{array}
$$
Here we use the usual notation for $p\in\R:$ $[p]=\max\,\{j\in\Z;\,\,j\leq p\}.$
\end{tm}

{\bf Proof:} $(a)$ These are all $\l\in\Lambda$ such that
$$
\l_1\in(m_1+k+\Z_+),\,\,\,\l_j=m_{j-1}+k-j+1\,\,\,2\leq j\leq k,\,\,\,\l_{k+1}=\pm m_k.
$$
\ind $(b)$ These are all $\l\in\Lambda$ such that
$$
\beg{array}{ll}
\l_s=m_s+k-s,&\,\,1\leq s\leq j-1,\\
\l_{j-1}>\l_j>\l_{j+1},&\\
\l_s=m_{s-1}+k-s+1,&\,\,j+1\leq s\leq k,\\
\l_{k+1}=\pm m_k.&
\end{array}
$$
\ind$(c)$ These are all $\l\in\Lambda$ such that
$$
\beg{array}{l}
\l_s=m_s+k-s,\,\,\,1\leq s\leq k,\,\,\,|\l_{k+1}|<m_k.
\end{array}
$$
The number of such $\l$'s is $0$ if $m_k=0$ or $m_k=\fr{1}{2},$ exactly $1$ if $m_k=1$ $(\l_{k+1}=0),$ and twice the number of natural numbers $<m_k$ if $m_k\geq\fr{3}{2}.$\\

\ind We now parametrize $\hat{G}^0.$ A class in $\wp{G}^0$ is unitary if and only if it is an irreducible subquotient of an end of complementary series. For $\s\in\hat{M}$ the complementary series is nonempty if and only if $\s$ is selfcontragredient, i.e. equivalent to its contragredient. Contragredient representation of\lb $\s=(n_1,\ldots,n_{k-1},n_k)$ is $-\s=(n_1,\ldots,n_{k-1},-n_k).$ Thus, $\s$ is selfcontragredient if and only if $n_k=0.$ In this case we set
$$
\nu(\s)=\min\,\{\nu\geq0;\,\,\pi^{\s,\nu}\,\,\mr{is}\,\,\mr{reducible}\}.
$$
From the necessary and sufficient conditions for reducibility of elementary representations we find that for $\s=(n_1,\ldots,n_{k-1},0)\in\hat{M}:$\\
\ind $\bullet$ If $n_1=\cdots=n_{k-1}=0,$ i.e. if $\s=\s_0=(0,\ldots,0)$ is the trivial onedimensional representation of $M,$ then
$$
\nu(\s_0)=k.
$$
In this case
$$
\G(\t^{\s_0,k})=\{(0,\ldots,0)\}\qu\mr{and}\qu\G(\o^{\s_0,k})=\{(s,0,\ldots,0);\,\,s\in\N\}
$$
and so $q^V(\t^{\s_0,k})=(0,\ldots,0)$ and $q^V(\o^{\s_0,k})=(1,0,\ldots,0).$\\
\ind $\bullet$ If $n_1>0,$ let $j\in\{2,\ldots,k\}$ be the smallest index such that $n_{j-1}>0.$ Then
$$
\nu(\s)=k-j+1.
$$
The $K-$spectra of irreducible subquotients of $\pi^{\s,k-j+1}$ consist of all $(m_1,\ldots,m_k)$ in $\hat{K}\cap\Z_+^k$ such that
$$
\beg{array}{ll}
\G(\t^{\s,k-j+1}):&\,\,m_1\geq n_1\geq\cdots\geq m_{j-1}\geq n_{j-1}\qu\mr{and}\qu m_s=0\,\,\fa s\geq j,\\
&\\
\G(\o^{\s,k-j+1}):&\,\,m_1\geq n_1\geq\cdots\geq m_{j-1}\geq n_{j-1}\geq m_j\geq1\qu\mr{and}\qu m_s=0\,\,\fa s>j.
\end{array}
$$
So we have
$$
q^V(\t^{\s,k-j+1})=(n_1,\ldots,n_{j-1},0,\ldots,0),\,\,\,q^V(\o^{\s,k-j+1})=(n_1,\ldots,n_{j-1},1,0,\ldots,0).
$$
Thus
\beg{tm} The map $\pi\map q^V(\pi)$ is a bijection of $\hat{G}^0$ onto
$$
\hat{K}_0=\{q=(m_1,\ldots,m_k)\in\hat{K};\,\,m_k=0\}.
$$
\end{tm}

\beg{thebibliography}{99}
\bibitem{pa} A. M. Gavrilik and A. U. Klimyk, \emph{Irreducible and indecomposable representations of the algebras $\ss\oo(n,1)$ and $\ii\ss\oo(n)$}, (in Russian) Preprint ITP$-$73$-$153R, Kiev, 1975.
\bibitem{pa} A. M. Gavrilik and A. U. Klimyk, \emph{Analysis of the representations of the Lorents and Euclidean groups of $n-$th order}, Preprint ITP$-$75$-$18E, Kiev, 1975.
\bibitem{pa} A. Guichardet, \emph{Repr\' esentations des groupes $\SO_0(n,1)$}, unpublished manuscript 1976.
\bibitem{pa} Harish$-$Chandra, \emph{On some applications of the universal enveloping algebra of a semi$-$simple Lie algebra}, Transactions of the American Mathematical Society, vol. 70(1951), pp. $28-96.$
\bibitem{pa} Harish$-$Chandra, \emph{Representations of a semi$-$simple Lie group on Banach spaces I}, Transactions of the American Mathematical Society, vol. 75(1953), pp. $185-243.$
\bibitem{pa} Harish$-$Chandra, \emph{Representations of a semi$-$simple Lie group on Banach spaces II}, Transactions of the American Mathematical Society, vol. 76(1954), pp. $26-65.$
\bibitem{pa} T. Hirai, \emph{On infinitesimal operators of irreducible representations of the Lorentz group of $n-$th order}, Proceedings of the Japan Academy, vol. 38(1962), pp. $83-87.$
\bibitem{pa} T. Hirai, \emph{On irreducible representations of the Lorentz group of $n-$th order}, Proceedings of the Japan Academy, vol. 38(1962), pp. $258-262.$
\bibitem{pa} A. U. Klimyk and A. M. Gavrilik, \emph{Representations matrix elements and Clebsch$-$Gordan coefficients of the semisimple Lie groups}, Journal of Mathematical Physics, vol. 20(1979), pp. $1624-1642.$
\bibitem{pa} H. Kraljevi\' c, \emph{On representations of the group $\SU(n,1)$}, Transactions of the American Mathematical Society, vol. 221(1976), pp. $433-448.$
\bibitem{pa} U. Ottoson, \emph{A classification of the unitary irreducible represnettions of $\SO_0(N,1)$}, Communications in Mathematical Physics, vol. 8(1968), pp. $228-244.$
\bibitem{pa} D. P. Zhelobenko, \emph{Description of quasisimple irreducible representations of the groups $\mr{U}(n,1),$ $\Spin(n,1)$} (in Russian), Izvestia Akademii nauk SSSR, Seria matematiceskaia, vol. 41(1977), pp. $34-53.$
\end{thebibliography}

\end{document}